\documentclass[12pt]{article}
\usepackage[T2A]{fontenc}
\usepackage[utf8]{inputenc}
\usepackage{amssymb,amsfonts,amsmath,color}
\usepackage{amssymb,amsmath,amsfonts,amsthm,amscd,latexsym,indentfirst,verbatim,mathtext,cite,enumerate,float}
\def\Span{\mathrm{Span}}
\def\Aut{\mathrm{Aut}}

\begin{document}

\sloppy

\begin{center}
{\Large
Unital decompositions of the matrix algebra of order three}

V. Gubarev
\end{center}

\begin{abstract}
We classify all decompositions of $M_3(\mathbb{C})$ into a direct vector-space sum of two
subalgebras such that one of the subalgebras contains the identity matrix.
\end{abstract}

\medskip
{\it Keywords}: decomposition of algebra, sum of rings, matrix algebra.

\section{Introduction}

Motivated by the results from group theory, O.H. Kegel proved in 1963~\cite{Kegel1963} that
an associative algebra represented as a vector-space sum of two nilpotent algebras is itself nilpotent.
After that, the study of decompositions of associative algebras into a sum of
two (not necessarily nilpotent) subalgebras began.
Let us mention the long-standing problem of K.I. Beidar and A.V. Mikhalev (\!\cite{BeidarMikhalev95},\,1995)
asking if a sum of two PI-algebras is again a~PI-algebra.
Recently, this problem was positively solved by M.~K\k{e}pczyk~\cite{Kepczyk2017}.
The famous K\"{o}the problem (If a ring $R$ has no nonzero nil ideals,
does it follow that $R$~has no nonzero nil one-sided ideals?)
is equivalent to a problem concerned decompositions~\cite{Ferrero1989}.

There is a series of works devoted to decompositions of Lie algebras~\cite{Goto,Onishchik69,Zusmanovich}.

It is naturally to study decompositions involved matrix algebras.
In~1999, Y.A.~Bahturin and O.H.~Kegel~\cite{BahturinKegel} described all algebras
decomposed as a sum of two matrix algebras. Moreover, they stated that
there are no proper decompositions of a~matrix algebra into a sum of two simple subalgebras.
In~\cite{BurdeGubarev}, it was stated that there are no proper semisimple
decompositions of the matrix algebra $M_n(\mathbb{C})$.

In~\cite{GonGub} (see also~\cite{AgoreCCP}), all direct decompositions of $M_2(\mathbb{C})$
were classified. The main aim of the current work is to get a classification
of all direct decompositions of $M_3(\mathbb{C})$ into a~sum of two subalgebras
such that one of the subalgebras contains the identity matrix.
Let us call such decomposition as a unital one.
Note that E.I. Konovalova in 2007 described~\cite{Konovalova}
all direct decompositions of $\mathrm{sl}_3(\mathbb{C})$ (25 cases).

Decompositions into a sum of two subalgebras are deeply
connected to Rota---Baxter operators~\cite{BurdeGubarev}.
In particular, given a direct decomposition of an algebra~$A$
one can construct a Rota---Baxter operator of nonzero weight on~$A$
(so called splitting Rota---Baxter operator).
In~\cite{GonGub}, all non-splitting Rota---Baxter operators
of nonzero weight on $M_3(\mathbb{C})$ were described (36 cases, all of them without parameters).
Thus, this work is the next step towards the complete classification
of Rota---Baxter operators of nonzero weight on $M_3(\mathbb{C})$.
After it, the final remaining step is to classify all direct decompositions of $M_3(\mathbb{C})$
into a sum of two subalgebras such that none of them contains the identity matrix.

We split the problem of classification of unital decompositions $M_3(\mathbb{C}) = A\oplus B$,
where $\dim B>\dim A$, into 3~main cases: $\dim B = 7$, $\dim B = 6$, or $\dim B = 5$.
In the first case, it is known that there is a~unique 7-dimensional subalgebra of $M_3(\mathbb{C})$
up to transpose and action of $\Aut(M_3(\mathbb{C}))$, it is a~subalgebra
$M = \Span\{e_{11},e_{12},e_{13},e_{22},e_{23},e_{32},e_{33}\}$~\cite{AgoreMax}.
Given a direct decomposition $M_3(\mathbb{C}) = S\oplus M$,
we may assume that~$S$ is generated by matrices $e_{21} + m_1$ and $e_{31} + m_2$,
where $m_1,m_2\in M$. Applying the condition that $S$ is a subalgebra, we find
possible decompositions.
To simplify such classification, it is very important to use automorphisms
of $M_3(\mathbb{C})$ preserving~$M$.

In~\cite{Unital}, it was proved that up to transpose and action of $\Aut(M_3(\mathbb{C}))$
there are two 6-dimensional subalgebras of $M_3(\mathbb{C})$:
either the subalgebra of upper-triangular matrices or the subalgebra of matrices with zero first column.
Finally, we state that there are exactly six different 5-dimensional subalgebras of
$M_3(\mathbb{C})$ (Lemma~5). However, both non-unital 5-dimensional subalgebras of $M_3(\mathbb{C})$
from Lemma~5 have the same collection of possible complement subalgebras~$S$.

So, in general we have eight subcases, and we get the complete classification
of such decompositions in Theorems~1--8.
Altogether there are 71 cases (when one separates decompositions with a unital
4-dimensional algebra from Theorem~4), some of them involve one or two parameters.
In some subcases (see Theorems~3,\,6, and 8) we are able to prove that
all obtained variants lie in different orbits under
automorphisms and antiautomorphisms (i.\,e., compositions of an automorphism and transpose)
of~$M_3(\mathbb{C})$.
Here we exclude the trivial decomposition $M_3(\mathbb{C}) = M_3(\mathbb{C}) \oplus (0)$.

\section{(7,2)-decompositions}

In what follows, we will apply an automorphism $\Theta_{12}$ of $M_3(\mathbb{C})$,
acting as follows, $\Theta_{12}(X) = T^{-1}XT$ for $T = e_{12}+e_{21}+e_{33}$.
Analogously, we define $\Theta_{13},\Theta_{23}\in\Aut(M_3(\mathbb{C}))$.

Every 2-dimensional associative algebra $L(e_1,e_2)$ over $\mathbb{C}$
is one of the following type~\cite{2Dim}, where we write down only nonzero products:

(D1) $N\oplus N$;

(D2) $N^2$: $e_1^2 = e_2$;

(D3) $F\oplus N$: $e_1^2 = e_1$;

(D4) $F\rtimes N_{unital}$: $e_1^2 = e_1$, $e_1e_2 = e_2e_1 = e_2$;

(D5) $F\rtimes N_{l,unital}$: $e_1^2 = e_1$, $e_1e_2 = e_2$;

(D6) $F\rtimes N_{r,unital}$: $e_1^2 = e_1$, $e_2e_1 = e_2$;

(D7) $F\oplus F$: $e_1^2 = e_1$, $e_2^2 = e_2$.

{\bf Theorem 1}.
Every direct decomposition of $M_3(\mathbb{C})$
with two subalgebras of the dimensions~2 and~7
up to transpose and up to action of $\Aut(M_3(\mathbb{C}))$ is isomorphic to
$S\oplus \Span\{e_{11},e_{12},e_{13},e_{22},e_{23},e_{32},e_{33}\}$,
where $S$ is one of the following subalgebras:

(R1) $\Span\{e_{21},e_{31}\}$;

(R2) $\Span\{e_{21},e_{31}+e_{23}\}$;

(R3) $\Span\{e_{21}+e_{22},e_{31}\}$;

(R4) $\Span\{e_{21}+e_{22}+e_{33},e_{31}+e_{32}\}$;

(R5) $\Span\{e_{21}+e_{22}+e_{33},e_{31}\}$;

(R6) $\Span\{e_{21}+e_{22},e_{31}+e_{32}\}$;

(R7) $\Span\{e_{21}+e_{22},e_{31}+e_{33}\}$;

(R8) $\Span\{e_{21}+(1-y)e_{22}+e_{23},e_{31}+e_{22}+e_{23}+ye_{33}\}$, $y\neq0$;

(R9) $\Span\{e_{21}+e_{22}+e_{33},e_{31}+e_{32}+e_{23}+ye_{33}\}$, $y\neq0$;

(R10) $\Span\left\{\begin{pmatrix}
0 & 0 & 0 \\
1 & d & 1 \\
0 & 1/f & 1 \\
\end{pmatrix},\begin{pmatrix}
0 & 0 & 0 \\
0 & 1 & f \\
1 & 1 & 1+f(1-d) \\
\end{pmatrix}\right\}$, $f\neq0$.

{\sc Proof}.
It is known~\cite{AgoreMax} that up to transpose and up to isomorphism there is the only 7-dimensional
subalgebra of $M_3(\mathbb{C})$, it is $M = \Span\{e_{11},e_{12},e_{13},e_{22},e_{23},e_{32},e_{33}\}$.
So, we want to describe all direct decompositions of $M_3(\mathbb{C})$
of the form $S\oplus M$, where $S$ is a~plane generated by matrices
\begin{equation}\label{v1v2}
v_1 = \begin{pmatrix}
a & b & c \\
1 & d & e \\
0 & f & g \\
\end{pmatrix},\quad
v_2 = \begin{pmatrix}
p & q & r \\
0 & s & t \\
1 & x & y \\
\end{pmatrix}.
\end{equation}

{\bf Lemma 1}.
Let $\varphi$ be an automorphism of $M_3(\mathbb{C})$ preserving the subalgebra~$M$.
Then $\varphi$ acts as follows,
\begin{gather}
e_{11}\to \begin{pmatrix}
1 & \beta & \gamma \\
0 & 0 & 0 \\
0 & 0 & 0 \\
\end{pmatrix},\quad
e_{12}\to \begin{pmatrix}
0 & \kappa & \lambda \\
0 & 0 & 0 \\
0 & 0 & 0 \\
\end{pmatrix},\quad
e_{13}\to \begin{pmatrix}
0 & \mu & \nu \\
0 & 0 & 0 \\
0 & 0 & 0 \\
\end{pmatrix}, \nonumber \\ \allowdisplaybreaks
e_{22}\to \frac{1}{\Delta}\begin{pmatrix}
0 & \kappa(\gamma \mu-\beta \nu) & \lambda(\gamma \mu-\beta \nu) \\
0 & \kappa \nu & \lambda \nu \\
0 & -\kappa \mu & -\lambda \mu \\
\end{pmatrix},\
e_{23}\to \frac{1}{\Delta}\begin{pmatrix}
0 & \mu(\gamma \mu-\beta \nu) & \nu(\gamma \mu-\beta \nu) \\
0 & \mu\nu & \nu^2 \\
0 & -\mu^2 & -\mu\nu \\
\end{pmatrix}, \nonumber \\
e_{32}\to \frac{1}{\Delta}\begin{pmatrix}
0 & \kappa(\beta \lambda-\gamma \kappa) & \lambda(\beta \lambda-\gamma \kappa) \\
0 & -\kappa\lambda & -\lambda^2 \\
0 & \kappa^2 & \kappa\lambda \\
\end{pmatrix},\
e_{33}\to \frac{1}{\Delta}\begin{pmatrix}
0 & \mu(\beta \lambda-\gamma \kappa) & \nu(\beta \lambda-\gamma \kappa) \\
0 & -\lambda\mu & -\lambda\nu \\
0 & \kappa\mu & \kappa\nu \\
\end{pmatrix}, \nonumber \\
e_{21}\to \frac{1}{\Delta}\begin{pmatrix}
\gamma \mu - \beta \nu & \beta(\gamma \mu -\beta \nu) & \gamma(\gamma \mu-\beta \nu) \\
\nu & \beta \nu & \gamma \nu \\
-\mu & -\beta \mu & -\gamma \mu \\
\end{pmatrix}, \label{AutPreserv} \\
e_{31}\to \frac{1}{\Delta}\begin{pmatrix}
\beta \lambda-\gamma \kappa & \beta(\beta \lambda-\gamma \kappa) & \gamma(\beta \lambda-\gamma \kappa) \\
-\lambda & -\beta \lambda & -\gamma \lambda \\
\kappa & \beta \kappa & \gamma \kappa \\
\end{pmatrix}, \nonumber
\end{gather}
where $\Delta = \kappa\nu-\lambda\mu\neq0$.

{\sc Proof}.
Since $\varphi$ has to preserve the radical of~$M$, we have
$e_{12}\to \kappa e_{12} + \lambda e_{13}$ and
$e_{13}\to \mu e_{12} + \nu e_{13}$ with $\Delta = \kappa\nu-\lambda\mu\neq0$.
Denote $\varphi(e_{11}) = X = (x_{ij})$.
Since $e_{11}$ acts as a left unit on $e_{12},e_{13}$,
we derive that $x_{11} = 1$.
On the other hand, $e_{1i}e_{11} = 0$ for $i=2,3$.
It implies that $x_{22} = x_{23} = x_{32} = x_{33} = 0$, and so,
$\varphi(e_{11}) = e_{11} + \beta e_{12} + \gamma e_{13}$.

Let $\varphi(e_{21}) = Z = (z_{ij})$.
From $e_{21}e_{11} = e_{21}$, we get that
$Z = \begin{pmatrix}
z_{11} & \beta z_{11} & \gamma z_{11} \\
z_{21} & \beta z_{21} & \gamma z_{21} \\
z_{31} & \beta z_{31} & \gamma z_{31} \\
\end{pmatrix}$.
The equality $e_{11}e_{21} = 0$ implies that
$z_{11} + \beta z_{21} + \gamma z_{31} = 0$.
From $e_{12}e_{21} = e_{11}$ and $e_{13}e_{21} = 0$,
we derive that
$\kappa z_{21} + \lambda z_{31} = 1$ and
$\mu z_{21} + \nu z_{31} = 0$. Solving this system,
we get $z_{21} = \nu/\Delta$, $z_{31} = -\mu/\Delta$, and
$z_{11} = (- \beta \nu + \gamma \mu)/\Delta$.

Analogously, we compute $\varphi(e_{31})$.
Then every element $\varphi(e_{ij})$ for $i,j\in\{2,3\}$
can be computed as $\varphi(e_{i1})\varphi(e_{1j})$.
\hfill $\square$

Let us show that we may apply
$\varphi(\beta,\gamma,\kappa,\lambda,\mu,\nu) = \varphi$~\eqref{AutPreserv}
such that $a = p = 0$. Indeed,
$$
\varphi(v_1) = \frac{1}{\Delta}\begin{pmatrix}
a\Delta + \gamma \mu - \beta \nu & * & * \\
\nu & * & * \\
-\mu & * & * \\
\end{pmatrix},\quad
\varphi(v_2) = \frac{1}{\Delta}\begin{pmatrix}
p\Delta + \beta \lambda - \gamma \kappa & * & * \\
-\lambda & * & * \\
\kappa & * & * \\
\end{pmatrix}.
$$
Then the following vectors generate the subalgebra $\varphi(\Span\{v_1,v_2\})$,
$$
\kappa\varphi(v_1) + \mu\varphi(v_2)
 = \begin{pmatrix}
\kappa a + \mu p - \beta & * & * \\
1 & * & * \\
0 & * & * \\
\end{pmatrix},\
\lambda\varphi(v_1) + \nu\varphi(v_2)
 = \begin{pmatrix}
\lambda a + \nu p - \gamma & * & * \\
0 & * & * \\
1 & * & * \\
\end{pmatrix}.
$$
Thus, it is enough to take $\beta = \kappa a + \mu p$ and $\gamma = \lambda a + \nu p$.
So, we may assume that $\underline{a = p = 0}$ in~\eqref{v1v2}.

Since $S$ is a subalgebra, we have the following multiplication table
$$
v_1^2 = dv_1 + fv_2,\quad
v_2^2 = tv_1 + yv_2,\quad
v_1v_2 = ev_1 + gv_2,\quad
v_2v_1 = sv_1 + xv_2.
$$
It gives rise to the equalities $b = c = q = r = 0$ and
the following system of equations
\begin{gather*}
\{f,t\}\times \{e-s,g-x\} = 0, \\
ft - eg = 0, \quad ft - sx = 0, \\
g(g-d) + f(e-y) = 0,\\
x(d-x) + f(y-s) = 0, \\ \allowdisplaybreaks
s(s-y) + t(x-d) = 0, \\
e(y-e) + t(d-g) = 0, \\
d(s-e) + ex - gs = 0, \quad
y(g-x) + ex - gs = 0.
\end{gather*}
Here by $X\times Y$ we mean a collection of the equations $x_iy_j = 0$,
where $X= \{x_i\}$ and $Y = \{y_j\}$.

{\sc Case I}: $f = t = 0$. Then we have the system
\begin{gather*}
eg = 0, \quad sx = 0, \\
g(d-g) = 0, \quad
x(d-x) = 0, \\
s(s-y) = 0, \quad
e(y-e) = 0, \\
d(s-e) + ex - gs = 0, \\
y(g-x) + ex - gs = 0.
\end{gather*}

{\sc Case IA}: $e = 0$.
So, the remaining system has the form
\begin{gather*}
sx = 0, \quad
\{g,s\}\times\{d-g\} = 0,\\
x(d-x) = 0, \quad
s(s-y) = 0, \quad
y(g-x) - gs = 0.
\end{gather*}

{\sc Case IAA}: $g = s = 0$.
Then $x(d-x) = 0$ and $xy = 0$.
If $x = 0$, then we get the first solution
$$
v_1 = \begin{pmatrix}
0 & 0 & 0 \\
1 & d & 0 \\
0 & 0 & 0 \\
\end{pmatrix},\quad
v_2 = \begin{pmatrix}
0 & 0 & 0 \\
0 & 0 & 0 \\
1 & 0 & y \\
\end{pmatrix}.
$$

If $x\neq0$, then $d = x$ and $y = 0$, it is the second solution
$$
v_1 = \begin{pmatrix}
0 & 0 & 0 \\
1 & x & 0 \\
0 & 0 & 0 \\
\end{pmatrix},\quad
v_2 = \begin{pmatrix}
0 & 0 & 0 \\
0 & 0 & 0 \\
1 & x & 0 \\
\end{pmatrix},\quad x\neq0.
$$

{\sc Case IAB}: $(g,s)\neq(0,0)$.
Then $d = g$ and
$$
sx = 0, \quad
x(g-x) = 0, \quad
s(s-y) = 0, \quad
g(s-y) + xy = 0.
$$

If $x = 0$, then $s = y$. If both $g,s$ are nonzero, we find the third solution
$$
v_1 = \begin{pmatrix}
0 & 0 & 0 \\
1 & d & 0 \\
0 & 0 & d \\
\end{pmatrix},\quad
v_2 = \begin{pmatrix}
0 & 0 & 0 \\
0 & y & 0 \\
1 & 0 & y \\
\end{pmatrix},\quad d,y\neq0.
$$
If one of $g,s$ is zero, up to the action of
$\Theta_{23} = \varphi(\beta=\gamma=\kappa=\nu=0,\lambda=\mu=1)$, we get the fourth solution
$$
v_1 = \begin{pmatrix}
0 & 0 & 0 \\
1 & d & 0 \\
0 & 0 & d \\
\end{pmatrix},\quad
v_2 = \begin{pmatrix}
0 & 0 & 0 \\
0 & 0 & 0 \\
1 & 0 & 0 \\
\end{pmatrix},\quad d\neq0.
$$

Now, consider when $x\neq0$. Then $s = 0$, $d = x\neq0$,
we write down the fifth solution,
$$
v_1 = \begin{pmatrix}
0 & 0 & 0 \\
1 & d & 0 \\
0 & 0 & d \\
\end{pmatrix},\quad
v_2 = \begin{pmatrix}
0 & 0 & 0 \\
0 & 0 & 0 \\
1 & d & y \\
\end{pmatrix},\quad d\neq0.
$$

{\sc Case IB}: $e \neq 0$.
Then $g = 0$, $y = e\neq0$, and
$$
sx = 0, \quad
x(d-x) = 0, \quad
s(s-e) = 0, \quad
d(s-e) + ex = 0.
$$
We may suppose that $x\neq0$, otherwise
the automorphism $\Theta_{23}$ converts this case to the case $e = 0$ and $x\neq0$.
So, $s = 0$ and $d = x$, it is the sixth solution
$$
v_1 = \begin{pmatrix}
0 & 0 & 0 \\
1 & x & e \\
0 & 0 & 0 \\
\end{pmatrix},\quad
v_2 = \begin{pmatrix}
0 & 0 & 0 \\
0 & 0 & 0 \\
1 & x & e \\
\end{pmatrix},\quad e,x\neq0.
$$

{\sc Case II}: $(f,t)\neq(0,0)$.
Then $e = s$ and $g = x$. Moreover,
$$
ft - eg = 0, \quad
g(d-g) + f(y-e) = 0, \quad
e(e-y) + t(g-d) = 0.
$$

{\sc Case IIA}: $f = 0$. Then
$$
eg = 0, \quad
g(d-g) = 0, \quad
e(e-y) + t(g-d) = 0.
$$
If $g = 0$, then $e(e-y) - td = 0$.
Since $t\neq0$, we get $d = e(e-y)/t$. It is the seventh solution
$$
v_1 = \begin{pmatrix}
0 & 0 & 0 \\
1 & \frac{e(e-y)}{t} & e \\
0 & 0 & 0 \\
\end{pmatrix},\quad
v_2 = \begin{pmatrix}
0 & 0 & 0 \\
0 & e & t \\
1 & 0 & y \\
\end{pmatrix},\quad t\neq0.
$$

If $g\neq0$, then we have $e = 0$ and $d = g$, it is the eighth solution
$$
v_1 = \begin{pmatrix}
0 & 0 & 0 \\
1 & g & 0 \\
0 & 0 & g \\
\end{pmatrix},\quad
v_2 = \begin{pmatrix}
0 & 0 & 0 \\
0 & 0 & t \\
1 & g & y \\
\end{pmatrix},\quad t,g\neq0.
$$

{\sc Case IIB}: $f \neq 0$.
If $t = 0$, then the automorphism $\Theta_{23}$
converts this case to the case $f = 0$ and $t\neq0$.
So, $t\neq0$ and we express $t = ex/f$ and
$y = e + x(x-d)/f$. It is the nineth solution
$$
v_1 = \left(\begin{matrix}
0 & 0 & 0 \\
1 & d & e \\
0 & f & x \\
\end{matrix}\right),\quad
v_2 = \left(\begin{matrix}
0 & 0 & 0 \\
0 & e & ex/f \\
1 & x & e + x(x-d)/f \\
\end{matrix}\right),\quad e,f,x\neq0.
$$

{\bf Lemma 2}.
Let $v_1,v_2$ be two matrices defined by~\eqref{v1v2} both with the first zero row.
For $\varphi$~\eqref{AutPreserv} taken with
$\beta = \gamma = \lambda = \mu = 0$, we have
$$
k\varphi(v_1) = \begin{pmatrix}
0 & 0 & 0 \\
1 & \kappa d & \nu e \\
0 & \kappa^2f/\nu & \kappa g \\
\end{pmatrix}, \quad
n\varphi(v_2) = \begin{pmatrix}
0 & 0 & 0 \\
0 & \nu s & t\nu^2/\kappa  \\
1 & \kappa x & \nu y \\
\end{pmatrix}.
$$

{\sc Proof}. It follows by the definition of~$\varphi$.
\hfill $\square$

Let us write down all nine cases.
\begin{align*}
& 1)\ v_1 = \begin{pmatrix}
0 & 0 & 0 \\
1 & d & 0 \\
0 & 0 & 0 \\
\end{pmatrix},\quad
&& v_2 = \begin{pmatrix}
0 & 0 & 0 \\
0 & 0 & 0 \\
1 & 0 & y \\
\end{pmatrix}.
\\
& 2)\ v_1 = \begin{pmatrix}
0 & 0 & 0 \\
1 & x & 0 \\
0 & 0 & 0 \\
\end{pmatrix},\quad
&& v_2 = \begin{pmatrix}
0 & 0 & 0 \\
0 & 0 & 0 \\
1 & x & 0 \\
\end{pmatrix},\quad x\neq0.
\\
& 3)\ v_1 = \begin{pmatrix}
0 & 0 & 0 \\
1 & d & 0 \\
0 & 0 & d \\
\end{pmatrix},\quad
&& v_2 = \begin{pmatrix}
0 & 0 & 0 \\
0 & y & 0 \\
1 & 0 & y \\
\end{pmatrix},\quad d,y\neq0.
\\
& 4)\ v_1 = \begin{pmatrix}
0 & 0 & 0 \\
1 & d & 0 \\
0 & 0 & d \\
\end{pmatrix},\quad
&& v_2 = \begin{pmatrix}
0 & 0 & 0 \\
0 & 0 & 0 \\
1 & 0 & 0 \\
\end{pmatrix},\quad d\neq0.
\\
& 5)\ v_1 = \begin{pmatrix}
0 & 0 & 0 \\
1 & x & 0 \\
0 & 0 & x \\
\end{pmatrix},\quad
&& v_2 = \begin{pmatrix}
0 & 0 & 0 \\
0 & 0 & 0 \\
1 & x & y \\
\end{pmatrix},\quad x\neq0.
\\
& 6)\ v_1 = \begin{pmatrix}
0 & 0 & 0 \\
1 & x & y \\
0 & 0 & 0 \\
\end{pmatrix},\quad
&& v_2 = \begin{pmatrix}
0 & 0 & 0 \\
0 & 0 & 0 \\
1 & x & y \\
\end{pmatrix},\quad x,y\neq0.
\\
& 7)\ v_1 = \begin{pmatrix}
0 & 0 & 0 \\
1 & \frac{e(e-y)}{t} & e \\
0 & 0 & 0 \\
\end{pmatrix},\quad
&& v_2 = \begin{pmatrix}
0 & 0 & 0 \\
0 & e & t \\
1 & 0 & y \\
\end{pmatrix},\quad t\neq0.
\end{align*}
\begin{align*}
& 8)\ v_1 = \begin{pmatrix}
0 & 0 & 0 \\
1 & x & 0 \\
0 & 0 & x \\
\end{pmatrix},\quad
&& v_2 = \begin{pmatrix}
0 & 0 & 0 \\
0 & 0 & t \\
1 & x & y \\
\end{pmatrix},\quad t,x\neq0.
\\
& 9)\
v_1 = \left(\begin{matrix}
0 & 0 & 0 \\
1 & d & e \\
0 & f & x \\
\end{matrix}\right),\quad
&&
v_2 = \left(\begin{matrix}
0 & 0 & 0 \\
0 & e & ex/f \\
1 & x & e + x(x-d)/f \\
\end{matrix}\right),\quad e,f,x\neq0.
\end{align*}

Now we apply Lemma~2 to all cases.

{\sc Case 1}.
Applying $\varphi$, we get
either $\Span\{e_{21},e_{31}\}\cong N\oplus N$ (R1),
or $\Span\{e_{21}+e_{22},e_{31}\}\cong F\oplus N$ ((R3), up to $\Theta_{23}$),
or $\Span\{e_{21}+e_{22},e_{31}+e_{33}\}\cong F^2$ (R7).

{\sc Case 2}.
Applying $\varphi$, we get $\Span\{e_{21}+e_{22},e_{31}+e_{32}\}\cong F\rtimes N_{r,unital}$ (R6).

{\sc Case 3}.
Applying $\varphi$, we get $\Span\{e_{21}+e_{22}+e_{33},e_{31}+e_{22}+e_{33}\}\cong F\rtimes N_{l,unital}$.
Up to action of~$\varphi$ defined by $\beta = \gamma = \mu = 0$ and $\kappa = \lambda = \nu = 1$, we get~(R5).

{\sc Case 4}.
Applying $\varphi$, we get
$\Span\{e_{21}+e_{22}+e_{33},e_{31}\}\cong F\rtimes N_{l,unital}$ (R5).

{\sc Case 5}.
Applying $\varphi$, we get
either $\Span\{e_{21}+e_{22}+e_{33},e_{31}+e_{32}\}\cong F\rtimes N_{unital}$ (R4)
or $\Span\{e_{21}+e_{22}+e_{33},e_{31}+e_{32}+e_{33}\}\cong F^2$.
Applying~$\varphi$ with $\beta = \gamma = \lambda = 0$ and $\kappa = \mu = \nu = 1$ to~(R7),
we get the last subcase.

{\sc Case 6}.
Applying $\varphi$, we get
$\Span\{e_{21}+e_{22}+e_{23},e_{31}+e_{32}+e_{33}\}\cong F\rtimes N_{r,unital}$.
The automorphism $\varphi$ defined by $\beta = \gamma = \mu = 0$ and $\kappa = \lambda = \nu = 1$
maps it to~(R6).

{\sc Case 7}.
Applying $\varphi$, we get when $e = 0$
either $\Span\{e_{21},e_{31}+e_{23}\}\cong N^2$ (R2)
or $\Span\{e_{21},e_{31}+e_{23}+e_{33}\}\cong F\oplus N$.
Applying~$\varphi$ with $\beta = \gamma = \lambda = 0$ and
$\kappa = \nu = -\mu = 1$ to~(R3), we obtain the last case.
When $e\neq0$, we have
$\Span\{e_{21}+(1-y)e_{22}+e_{23},e_{31}+e_{22}+e_{23}+ye_{33}\}$~(R8).

{\sc Case 8}.
Applying $\varphi$, we get the subalgebra
$S(y) = \Span\{e_{21}+e_{22}+e_{33},e_{31}+e_{32}+e_{23}+ye_{33}\}$.
If $y = 0$, we have $\Span\{e_{21}+e_{22}+e_{33},e_{31}+e_{32}+e_{23}\}\cong F^2$.
Applying~$\varphi$ with $\beta = \gamma = \mu = 0$, $\kappa = \lambda = 1$, and $\nu = -2$
to~(R7), we get the case. Otherwise, it is~(R9).

{\sc Case 9}.
Applying $\varphi$, we get (R10).

Theorem is proved.\hfill $\square$

{\bf Remark 1}.
All cases (R1)--(R7) from Theorem~1 lie in different orbits
under action of automorphisms or antiautomorphisms of $M_3(\mathbb{C})$ preserving~$M$.
Indeed, a subalgebra~$S$ from case~($Ri$), where $i = 1,\ldots,7$,
is isomorphic to the two-dimensional algebra of type~($Di$).
Cases~(R5) and~(R6) are not still antiisomorphic, since
their idempotents have different ranks.

\section{Unital (6,3)-decompositions}

Let $A$ be a six-dimensional subalgebra of $M_3(\mathbb{C})$.
In~\cite{Unital}, it was actually proved that if $A$ is unital,
then up to transpose $A$ is isomorphic to the subalgebra of upper-triangular matrices.
Suppose that $A$ is not unital. Then $A$ consists only of degenerate matrices,
and then up to transpose~$A$ is isomorphic to the subalgebra having all zero elements in the first column~\cite{Meshulam}.
We will study both cases and describe decompositions with corresponding six-dimensional subalgebra.

\subsection{Case of unital 3-dimensional subalgebra}

Thus, we fix six-dimensional subalgebra $M = \Span\{e_{12},e_{13},e_{22},e_{23},e_{32},e_{33}\}$.

{\bf Theorem 2}.
Every direct decomposition of $M_3(\mathbb{C})$
with two subalgebras of the dimensions~3 and~6,
where a 3-dimensional subalgebra is unital, up to transpose
and up to action of $\Aut(M_3(\mathbb{C}))$ is isomorphic to
$S\oplus \Span\{e_{12},e_{13},e_{22},e_{23},e_{32},e_{33}\}$,
where $S$ is one of the following subalgebras:

(S1) $\Span\{E,e_{21},e_{31}\}$;

(S2) $\Span\{E,e_{21},e_{31}+e_{23}\}$;

(S3) $\Span\{E,e_{21}+e_{22},e_{31}\}$;

(S4) $\Span\{E,e_{21}+e_{22}+e_{33},e_{31}+e_{32}\}$;

(S5) $\Span\{E,e_{21}+e_{22}+e_{33},e_{31}\}$;

(S6) $\Span\{E,e_{21}+e_{22},e_{31}+e_{32}\}$;

(S7) $\Span\{E,e_{21}+e_{22},e_{31}+e_{33}\}$;

(S8) $\Span\{E,e_{21}+(1-y)e_{22}+e_{23},e_{31}+e_{22}+e_{23}+ye_{33}\}$, $y\neq0$;

(S9) $\Span\{E,e_{21}+e_{22}+e_{33},e_{31}+e_{32}+e_{23}+ye_{33}\}$, $y\neq0$;

(S10) $\Span\left\{E,\begin{pmatrix}
0 & 0 & 0 \\
1 & d & 1 \\
0 & 1/f & 1 \\
\end{pmatrix},\begin{pmatrix}
0 & 0 & 0 \\
0 & 1 & f \\
1 & 1 & 1+f(1-d) \\
\end{pmatrix}\right\}$, $f\neq0$;

(S11) $\Span\{E,e_{21}+e_{13}+de_{23}+e_{32},e_{31}+e_{12}+e_{23}+d(e_{22}+e_{33})\}$;

(S12) $\Span\left\{E,
\begin{pmatrix}
0 & 0 & eu-1 \\
1 & 1 & 1 \\
0 & e & 1 \\
\end{pmatrix},
\begin{pmatrix}
0 & eu-1 & 0 \\
0 & 1 & u \\
1 & 1 & 1 \\
\end{pmatrix}\right\}$, $e\neq0$, $eu\neq1$.

{\sc Proof}. A subalgebra~$S$ is generated by~$E$ and by matrices
\begin{equation}\label{v1v2:6-3}
v_1 = \begin{pmatrix}
0 & a & b \\
1 & c & d \\
0 & e & f \\
\end{pmatrix},\quad
v_2 = \begin{pmatrix}
0 & r & s \\
0 & t & u \\
1 & x & y \\
\end{pmatrix}.
\end{equation}

It is clear that given a decomposition $S'\oplus M'$ from Theorem~1,
one get a decomposition $S\oplus M$ with $S = S'\oplus\Span\{E\}$.
Thus, we have automatically decompositions~(S1)--(S10).

Let us note that if $\chi$ is an automorphism of $M_3(\mathbb{C})$
preserving~$M$, then $\chi$ has to preserve the subalgebra
$\Span\{e_{11},e_{12},e_{13},e_{22},e_{23},e_{32},e_{33}\}$.
Indeed, it follows from the equality $U_3 = \Span\{E\}\oplus M$.
All such automorphisms are described in Lemma~1.

{\bf Lemma 3}.
Up to action of an appropriate~$\varphi$~\eqref{AutPreserv}
we have either $a = s = 0$ or $a = 0$, $b+r = 0$.

{\sc Proof}.
Applying $\varphi$~\eqref{AutPreserv} with $\beta = \gamma = 0$, we get
\begin{gather*}
\tilde{v}_1 = \kappa v_1 + \mu v_2
= \begin{pmatrix}
0 & \kappa^2 a + \kappa \mu(b+r) + \mu^2 s & \kappa \lambda a + \kappa \nu b + \lambda \mu r + \mu \nu s \\
1 & * & * \\
0 & * & * \\
\end{pmatrix}, \\
\tilde{v}_2 = \lambda v_1 + \nu v_2 = \begin{pmatrix}
0 & \kappa\lambda a + \lambda \mu b + \kappa \nu r + \mu \nu s & \lambda^2a + \lambda\nu(b+r) + \nu^2s \\
0 & * & * \\
1 & * & * \\
\end{pmatrix}.
\end{gather*}
Since we always may choose $\kappa,\lambda,\mu,\nu$ such that $\Delta \neq 0$ and
$\kappa^2a + \kappa \mu(b+r) + \mu^2 s = 0$, we may assume that $a = 0$.

If $s = 0$ or $b + r = 0$, then we are done.
Suppose that $s\neq0$ and $b+r\neq0$, then
we may choose nonzero $\kappa$ and $\mu$ such that
$\kappa(b+r) + \mu s = 0$ and we take $\nu = 0$.
Thus, $\tilde{v}_1$ and $\tilde{v}_2$ have zero elements $a,s$.
\hfill $\square$

By Lemma~3, we assume that $\underline{a = 0}$ in~\eqref{v1v2:6-3}.
Further, we will consider different cases when $s = 0$ or $b+r = 0$.

Since $S$ is a subalgebra, we have the following multiplication table
$$
v_1^2 = cv_1 + ev_2,\quad
v_2^2 = uv_1 + yv_2 + sE,\quad
v_1v_2 = dv_1 + fv_2 + bE,\quad
v_2v_1 = tv_1 + xv_2 + rE.
$$
It gives rise to the following system of equations
\begin{gather}
\{e,u\}\times \{b-r,d-t,f-x\} = 0, \quad
bx - fr = 0, \nonumber \\
b + df - eu = 0, \quad
r + tx - eu = 0, \nonumber \\
bc - bf + es = 0, \quad
rc - rx + es  = 0, \nonumber \\
bd - by + sf = 0, \nonumber \quad
rt - ry + sx = 0, \label{6-3-I-system} \\
f^2 - cf + de - ey = 0, \nonumber \quad
x^2 - cx + te - ey = 0, \nonumber \\
t^2 - ty + ux - uc - s = 0, \quad
d^2 - dy + uf - uc - s = 0, \nonumber \\
dr + sf - sx - bt = 0, \nonumber \\
r + dx + ct - b - ft - cd = 0, \nonumber \\
b + dx + fy - r - ft - xy = 0. \nonumber
\end{gather}

{\sc Case I}: $s = 0$.

{\sc Case IA}: $e = u  = 0$.
Then the system~\eqref{6-3-I-system} transforms to the following one,
\begin{equation}\label{6-3-I-system1A}
\begin{gathered}
r(c - x)  = 0, \quad
r(t - y) = 0, \\
x(c - x) = 0, \quad
t(t - y) = 0, \\
b(d - y) = 0, \quad
b(c - f) = 0, \\
d(d - y)  = 0, \quad
f(c - f) = 0, \\
b + df = 0, \quad
r + tx = 0, \\
bx - fr = 0, \quad
dr - bt = 0, \\
r + dx + ct - b - ft - cd = 0, \\
b + dx + fy - r - ft - xy = 0.
\end{gathered}
\end{equation}

If $b = r = 0$, then we have exactly the same cases of~$S$ which are listed in
Theorem~1 joint with the identity matrix.
The cases when $b = 0$, $r\neq0$ and $b\neq0$, $r=0$ are mapped to each other by the automorphism~$\Theta_{23}$.

{\sc Case IAA}: $b\neq0$ and $r=0$.
By~\eqref{6-3-I-system1A} we have $d = y$, $c = f$, $x = t = 0$ and $b = - fy$, it's the first solution
$$
v_1 = \begin{pmatrix}
0 & 0 & -fy \\
1 & f & y \\
0 & 0 & f \\
\end{pmatrix},\quad
v_2 = \begin{pmatrix}
0 & 0 & 0 \\
0 & 0 & 0 \\
1 & 0 & y \\
\end{pmatrix},\quad f,y\neq0.
$$

{\sc Case IAB}: $b,r\neq0$.
By~\eqref{6-3-I-system1A} we have $c = f = x$, $d = y = t$.
Moreover, $b = r = - fy$ , it's the second solution
$$
v_1 = \begin{pmatrix}
0 & 0 & -fy \\
1 & f & y \\
0 & 0 & f \\
\end{pmatrix},\quad
v_2 = \begin{pmatrix}
0 & -fy & 0 \\
0 & y & 0 \\
1 & f & y \\
\end{pmatrix},\quad f,y\neq0.
$$

{\sc Case IB}: $(e,u)\neq(0,0)$.
Then we have $b = r$, $d = t$, $f = x$.
As above, we may assume that $b\neq0$.
So, $d = y$ and $c = x$.
Thus, the system~\eqref{6-3-I-system} converts to the only equality
$b + dx - eu = 0$. It is the third solution
$$
v_1 = \begin{pmatrix}
0 & 0 & eu-dx \\
1 & x & d \\
0 & e & x \\
\end{pmatrix},\quad
v_2 = \begin{pmatrix}
0 & eu-dx & 0 \\
0 & d & u \\
1 & x & d \\
\end{pmatrix},\quad (e,u)\neq(0,0),\ eu-dx\neq0.
$$

{\sc Case II}: $b + r = 0$. Not to repeat Case 1, we may assume that $s\neq0$.

{\sc Case IIA}: $e = u = 0$. Then the system~\eqref{6-3-I-system} has the form
\begin{gather*}
b(c - f) = 0, \quad
b(c - x) = 0, \quad
b(x + f) = 0, \\
f(c - f) = 0, \quad
x(x - c) = 0, \\
b + df = 0, \quad
b - tx = 0, \\
bd - by + sf = 0, \quad
by - bt + sx = 0, \\ \allowdisplaybreaks
t^2 - ty - s = 0, \quad
d^2 - dy - s = 0, \\
-bd + sf - sx - bt = 0, \\
-2b + dx + ct - ft - cd = 0, \\
2b + dx + fy - ft - xy = 0.
\end{gather*}

{\sc Case IIAA}: $b\neq0$.
Then $c = f = x = 0$ and $b = - df = 0$, a contradiction.

{\sc Case IIAB}: $b = 0$. Then $f = x = 0$, and
$$
c(t - d) = 0, \quad
t^2 - ty - s = 0, \quad
d^2 - dy - s = 0.
$$

If $d = t$, then it is the fourth solution
$$
v_1 = \begin{pmatrix}
0 & 0 & 0 \\
1 & c & d \\
0 & 0 & 0 \\
\end{pmatrix},\quad
v_2 = \begin{pmatrix}
0 & 0 & d(d-y) \\
0 & d & 0 \\
1 & 0 & y \\
\end{pmatrix},\quad d\neq0,\ d\neq y.
$$

If $d \neq t$, then it is the fifth solution
$$
v_1 = \begin{pmatrix}
0 & 0 & 0 \\
1 & 0 & d \\
0 & 0 & 0 \\
\end{pmatrix},\quad
v_2 = \begin{pmatrix}
0 & 0 & -dt \\
0 & t & 0 \\
1 & 0 & d+t \\
\end{pmatrix},\quad d,t\neq0.
$$

{\sc Case IIB}: $(e,u)\neq(0,0)$.
Then $b = r = 0$, $d = t$, $f = x = 0$, $e = 0$, $u\neq0$, since $s\neq0$.
We have the only remaining equality
$d^2 - dy - uc - s = 0$.
It is the sixth solution
$$
v_1 = \begin{pmatrix}
0 & 0 & 0 \\
1 & c & d \\
0 & 0 & 0 \\
\end{pmatrix},\quad
v_2 = \begin{pmatrix}
0 & 0 & d^2-dy-cu \\
0 & d & u \\
1 & 0 & y \\
\end{pmatrix},\quad u\neq0,\ d^2-dy-cu\neq0.
$$

Let us gather all obtained cases except ones arisen from $(7,2)$-decompositions:
\begin{align*}
& 1)\
v_1 = \begin{pmatrix}
0 & 0 & -fy \\
1 & f & y \\
0 & 0 & f \\
\end{pmatrix},\quad
&& v_2 = \begin{pmatrix}
0 & 0 & 0 \\
0 & 0 & 0 \\
1 & 0 & y \\
\end{pmatrix},\quad f,y\neq0.
\\
& 2)\
v_1 = \begin{pmatrix}
0 & 0 & -fy \\
1 & f & y \\
0 & 0 & f \\
\end{pmatrix},\quad
&& v_2 = \begin{pmatrix}
0 & -fy & 0 \\
0 & y & 0 \\
1 & f & y \\
\end{pmatrix},\quad f,y\neq0.
\\
& 3)\ v_1 = \begin{pmatrix}
0 & 0 & eu-dx \\
1 & x & d \\
0 & e & x \\
\end{pmatrix},\quad
&& v_2 = \begin{pmatrix}
0 & eu-dx & 0 \\
0 & d & u \\
1 & x & d \\
\end{pmatrix},\quad (e,u)\neq(0,0),\ eu-dx\neq0.
\\
& 4)\
v_1 = \begin{pmatrix}
0 & 0 & 0 \\
1 & c & d \\
0 & 0 & 0 \\
\end{pmatrix},\quad
&& v_2 = \begin{pmatrix}
0 & 0 & d(d-y) \\
0 & d & 0 \\
1 & 0 & y \\
\end{pmatrix},\quad d\neq0,\ d\neq y.
\\
& 5)\
v_1 = \begin{pmatrix}
0 & 0 & 0 \\
1 & 0 & d \\
0 & 0 & 0 \\
\end{pmatrix},\quad
&& v_2 = \begin{pmatrix}
0 & 0 & -dt \\
0 & t & 0 \\
1 & 0 & d+t \\
\end{pmatrix},\quad d,t\neq0.
\\
& 6)\
v_1 = \begin{pmatrix}
0 & 0 & 0 \\
1 & c & d \\
0 & 0 & 0 \\
\end{pmatrix},\quad
&& v_2 = \begin{pmatrix}
0 & 0 & d^2-dy-cu \\
0 & d & u \\
1 & 0 & y \\
\end{pmatrix},\quad u\neq0,\ d^2-dy-cu\neq0.
\end{align*}

{\bf Lemma 4}.
Let $v_1,v_2$ be two matrices defined by~\eqref{v1v2:6-3}.
For $\varphi$~\eqref{AutPreserv} taken with $\beta = \gamma = \lambda = \mu = 0$, we have
$$
\kappa\varphi(v_1) = \begin{pmatrix}
0 & \kappa^2 a & \kappa\nu b \\
1 & \kappa c & \nu d \\
0 & \kappa^2 e/\nu & \kappa f \\
\end{pmatrix}, \quad
n\varphi(v_2) = \begin{pmatrix}
0 & \kappa\nu r & \nu^2 s \\
0 & \nu t & \nu^2 u/\kappa \\
1 & \kappa x & \nu y \\
\end{pmatrix}.
$$

{\sc Proof}. It follows by the definition of~$\varphi$.
\hfill $\square$

{\sc Case 1}.
Applying $\varphi$ with $\beta = \gamma = \lambda = \mu = 0$,
$\kappa=1/f$, $\nu=1/y$, we get
$$
v_1 = \begin{pmatrix}
0 & 0 & -1 \\
1 & 1 & 1 \\
0 & 0 & 1 \\
\end{pmatrix},\quad
v_2 = \begin{pmatrix}
0 & 0 & 0 \\
0 & 0 & 0 \\
1 & 0 & 1 \\
\end{pmatrix}.
$$
Now, applying~$\varphi$ with
$\kappa = \mu = \nu = 1$, $\beta = \gamma = -1$, and $\lambda = 0$, we get
$\varphi(S) = \Span\{E,e_{21}-e_{31},e_{31}-e_{22}-e_{33}\}$,
which is one of the cases from Theorem~1.

{\sc Case 2}.
Applying $\varphi$ with
$\beta = \gamma = \lambda = \mu = 0$, $\kappa=1/f$, $\nu=1/y$, we get
$$
v_1 = \begin{pmatrix}
0 & 0 & -1 \\
1 & 1 & 1 \\
0 & 0 & 1 \\
\end{pmatrix},\quad
v_2 = \begin{pmatrix}
0 & -1 & 0 \\
0 & 1 & 0 \\
1 & 1 & 1 \\
\end{pmatrix}.
$$
Now, applying~$\varphi$ with
$\kappa = \mu = \nu = 1$, $\beta = \gamma = -1$, and $\lambda = 0$, we get
$\varphi(S) = \Span\{E,e_{21}-e_{31},e_{31}-e_{22}-2e_{33}\}$,
which is one of the cases from Theorem~1.

{\sc Case 3}.
If $x,d\neq0$, then we apply $\varphi$ with
$\beta = \gamma = \lambda = \mu = 0$, $\kappa = 1/x$, $\nu = 1/d$ and we get
$$
v_1 = \begin{pmatrix}
0 & 0 & eu-1 \\
1 & 1 & 1 \\
0 & e & 1 \\
\end{pmatrix},\quad
v_2 = \begin{pmatrix}
0 & eu-1 & 0 \\
0 & 1 & u \\
1 & 1 & 1 \\
\end{pmatrix},\quad (e,u)\neq(0,0),\ eu\neq1.
$$
Further, if $e = 0$, we apply~$\varphi$ with $\beta + \mu = \gamma + \nu = 0$
to get vectors with zero $e_{12}$- and $e_{13}$-coordinates,
this case is up to action of $\Aut(M_3(\mathbb{C}))$ one of the cases from Theorem~1.
Otherwise, it is~(S12).

If $x = 0$, then $eu\neq0$ and we apply
$\varphi$ with $\beta = \gamma = \lambda = \mu = 0$, $\kappa^2 e/\nu=1$, $\nu^2 u/\kappa=1$,
$$
v_1 = \begin{pmatrix}
0 & 0 & 1 \\
1 & 0 & d \\
0 & 1 & 0 \\
\end{pmatrix},\quad
v_2 = \begin{pmatrix}
0 & 1 & 0 \\
0 & d & 1 \\
1 & 0 & d \\
\end{pmatrix}.
$$
It is~(S11).

If $d = 0$, then up to the action of $\Theta_{23}$, we get the same subalgebra.

{\sc Cases 4~and~5}.
By~Lemma~4, we may assume that $d = 1$.
Now, applying~$\varphi$ with $\beta + \mu = 0$ and $\gamma + \nu = 0$, we get vectors
with zero $e_{12}$- and $e_{13}$-coordinates,
so, we get one of the cases from Theorem~1.

{\sc Case 6}.
If $c = 0$, then $d\neq0$, and we apply $\varphi$ with $\beta = \gamma = \lambda = \mu = 0$,
$\nu = 1/d$, $\kappa = u/d^2$ and we get
$$
v_1 = \begin{pmatrix}
0 & 0 & 0 \\
1 & 0 & 1 \\
0 & 0 & 0 \\
\end{pmatrix},\quad
v_2 = \begin{pmatrix}
0 & 0 & 1-y \\
0 & 1 & 1 \\
1 & 0 & y \\
\end{pmatrix},\quad y\neq1.
$$
Now, applying~$\varphi$ with $\beta + \mu = 0$ and $\gamma + \nu = 0$, we get vectors
with zero $e_{12}$- and $e_{13}$-coordinates,
it is one of the cases from Theorem~1.

If $c\neq0$, by Lemma~4, we may assume that $c = u = 1$.
Applying~$\varphi$ with $\beta = -\mu d - \kappa$ and $\gamma = - \nu d - \lambda$,
we get vectors with zero $e_{12}$- and $e_{13}$-coordinates,
it is one of the cases from Theorem~1.

Theorem is proved.\hfill $\square$

{\bf Remark 2}.
The decompositions $M_3(\mathbb{C}) = S\oplus \Span\{e_{12},e_{13},e_{22},e_{23},e_{32},e_{33}\}$
of ty\-pes~(S1) and~(S2) have actually appeared in~\cite[Example 2.7.4]{AgoreCCP}.

\subsection{Case of unital 6-dimensional subalgebra}

Thus, we fix six-dimensional subalgebra $M = \Span\{e_{11},e_{12},e_{13},e_{22},e_{23},e_{33}\}$.
Now we want to describe all direct decompositions of $M_3(\mathbb{C})$
of the form $S\oplus M$, where $S$ is a~3-dimensional subalgebra.

{\bf Theorem 3}.
Every direct decomposition of $M_3(\mathbb{C})$
with two subalgebras of the dimensions~3 and~6,
where a 6-dimensional subalgebra is unital,
up to transpose and up to action of $\Aut(M_3(\mathbb{C}))$
is isomorphic to
$S\oplus \Span\{e_{11},e_{12},e_{13},e_{22},e_{23},e_{33}\}$,
where $S$ is one of the following subalgebras:

(T1) $\Span\{e_{21},e_{31},e_{32}\}$;

(T2) $\Span\{e_{21}+e_{22},e_{31},e_{32}\}$;

(T3) $\Span\{e_{21}+e_{22}+e_{33},e_{31},e_{32}\}$;

(T4) $\Span\{e_{21}+e_{22},e_{11}+e_{22}+e_{31},e_{31}+e_{32}\}$;

(T5) $\Span\{e_{21}+e_{22},e_{11}+e_{22}+e_{31},e_{12}+e_{21}+e_{32}\}$;

(T6) $\Span\{e_{21}+e_{22}+e_{13}-e_{23}+e_{33},e_{11}+e_{22}+e_{31},e_{12}-e_{22}+e_{32}\}$.

{\sc Proof}. The subalgebra~$S$ is generated by matrices
$$
v_1 = \begin{pmatrix}
a & b & c \\
1 & d & e \\
0 & 0 & f \\
\end{pmatrix},\quad
v_2 = \begin{pmatrix}
g & h & i \\
0 & j & k \\
1 & 0 & l \\
\end{pmatrix},\quad
v_3 = \begin{pmatrix}
m & n & s \\
0 & p & q \\
0 & 1 & r \\
\end{pmatrix}.
$$

It is well-known that all automorphisms of $M$ are inner, so
$\psi\in\Aut(M_3(\mathbb{C}))$ preserving~$M$ has the form
\begin{equation}\label{AutU}
\begin{gathered}
e_{11}\to \begin{pmatrix}
1 & \beta & \gamma \\
0 & 0 & 0 \\
0 & 0 & 0 \\
\end{pmatrix},\quad
e_{12}\to \begin{pmatrix}
0 & \delta & \varepsilon \\
0 & 0 & 0 \\
0 & 0 & 0 \\
\end{pmatrix},\quad
e_{13}\to \begin{pmatrix}
0 & 0 & \alpha \\
0 & 0 & 0 \\
0 & 0 & 0 \\
\end{pmatrix}, \\
e_{21}\to
\frac{1}{\delta}\begin{pmatrix}
-\beta & -\beta^2 & -\beta\gamma \\
1 & \beta & \gamma \\
0 & 0 & 0 \\
\end{pmatrix},\quad
e_{22}\to
\begin{pmatrix}
0 & -\beta & -\beta\varepsilon/\delta \\
0 & 1 & \varepsilon/\delta \\
0 & 0 & 0 \\
\end{pmatrix}, \\
e_{23}\to
\frac{1}{\delta}\begin{pmatrix}
0 & 0 & -\alpha\beta \\
0 & 0 & \alpha \\
0 & 0 & 0 \\
\end{pmatrix},\quad
e_{31}\to
\frac{1}{\alpha\delta}
\begin{pmatrix}
\beta\varepsilon-\gamma\delta & \beta(\beta\varepsilon-\gamma\delta)
 & \gamma(\beta\varepsilon-\gamma\delta) \\
-\varepsilon & -\beta\varepsilon & -\gamma\varepsilon \\
\delta & \beta\delta & \gamma\delta \\
\end{pmatrix}, \\
e_{32}\to
\frac{1}{\alpha}\begin{pmatrix}
0 & \beta\varepsilon-\gamma\delta & \varepsilon(\beta\varepsilon-\gamma\delta)/\delta \\
0 & -\varepsilon & -\varepsilon^2/\delta \\
0 & \delta & \varepsilon \\
\end{pmatrix}, \quad
e_{33}\to
\frac{1}{\delta}\begin{pmatrix}
0 & 0 & \beta\varepsilon-\gamma\delta \\
0 & 0 & -\varepsilon \\
0 & 0 & \delta \\
\end{pmatrix},
\end{gathered}
\end{equation}
where $\alpha,\delta \neq0$.
Note that it is exactly~$\varphi$~\eqref{AutPreserv} considered with $\mu = 0$.

Let us act by $\psi$~\eqref{AutU} on the vectors~$v_1,v_2,v_3$:
\begin{gather*}
\psi(v_1)
= \begin{pmatrix}
a-\beta/\delta & * & * \\
1/\delta & * & * \\
0 & 0 & f
\end{pmatrix},\quad
\psi(v_3)
= \begin{pmatrix}
m & * & * \\
0 & * & * \\
0 & \delta/\alpha & r+\varepsilon/\alpha
\end{pmatrix},\\
\psi(v_2)
= \begin{pmatrix}
* & * & * \\
-\frac{\varepsilon}{\alpha\delta} & * & * \\
\frac{1}{\alpha} & \frac{\beta}{\alpha} & l+\frac{\gamma}{\alpha} \\
\end{pmatrix}
 = -\frac{\varepsilon v_1}{\alpha\delta}
 + \frac{\beta v_3}{\alpha}
 + \begin{pmatrix}
* & * & * \\
0 & * & * \\
\frac{1}{\alpha} & 0 & l+\frac{\gamma}{\alpha}
 + \frac{\varepsilon f}{\alpha\delta} - \frac{\beta r}{\alpha}
\end{pmatrix}.
\end{gather*}

Taking $\beta = \delta a$, $\varepsilon = -\alpha r$,
and $\gamma = \beta r - \varepsilon f/\delta - \alpha l$,
we may assume that $\underline{a = l = r = 0}$ up to the action of such $\psi$.

Now, we apply that $S$ is a subalgebra, so it is closed under the mutliplication. From
$$
v_1^2 = dv_1, \quad
v_2^2 = kv_1 + gv_2 + hv_3, \quad
v_3^2 = p v_3,
$$
we get $b = q = 0$ and the system
\begin{gather*}
f(d-f) = 0, \quad
c(d-f) = 0, \quad
ef + c = 0, \\
i = hm = kf,\quad
h(j-n) = 0, \quad
kh - ck - hs = 0, \\ \allowdisplaybreaks
j^2-jg-dk-hp = 0, \quad
k(j-e-g) = 0, \\
m(m-p) = 0, \quad
s(m-p) = 0, \quad
mn + s = 0.
\end{gather*}

Currently, we have
$$
v_1 = \begin{pmatrix}
0 & 0 & -ef \\
1 & d & e \\
0 & 0 & f \\
\end{pmatrix},\quad
v_2 = \begin{pmatrix}
g & h & kf \\
0 & j & k \\
1 & 0 & 0 \\
\end{pmatrix},\quad
v_3 = \begin{pmatrix}
m & n & -mn \\
0 & p & 0 \\
0 & 1 & 0 \\
\end{pmatrix},
$$
and the system
\begin{gather*}
f(d-f) = 0, \quad
kf - hm = 0,\quad
h(j-n) = 0, \quad
j^2-jg-dk-hp = 0, \\
k(h + ef + fn) = 0, \quad
k(j-e-g) = 0,\quad
m(m-p) = 0.
\end{gather*}

From
$v_1v_3 = mv_1+fv_3$ and
$v_3v_1 = pv_1+v_2+dv_3$, we obtain $h = k = 0$, and
\begin{gather*}
mf = 0, \quad
f(e+n) = 0, \quad
m(e+n) = 0, \quad
n+e+p(d-f) - md = 0, \\
g = n - md, \quad
j = -pd, \quad
e = pf, \quad
en + mnd + epf = 0.
\end{gather*}

We rewrite $v_1,v_2,v_3$ as follows,
$$
v_1 = \begin{pmatrix}
0 & 0 & -pf^2 \\
1 & d & pf \\
0 & 0 & f \\
\end{pmatrix},\quad
v_2 = \begin{pmatrix}
-pd & 0 & 0 \\
0 & -pd & 0 \\
1 & 0 & 0 \\
\end{pmatrix},\quad
v_3 = \begin{pmatrix}
m & d(m-p) & 0 \\
0 & p & 0 \\
0 & 1 & 0 \\
\end{pmatrix}.
$$
Moreover, we have the system
$$
mf = 0, \quad
f(d-f) = 0, \quad
m(m-p) = 0.
$$

The products between $v_2$ and $v_i$, $i=1,3$,
give us no new equations.

{\sc Case I}. $f = 0$.

{\sc Case IA}. $m = 0$.
We have the first solution
$$
v_1 = \begin{pmatrix}
0 & 0 & 0 \\
1 & d & 0 \\
0 & 0 & 0 \\
\end{pmatrix},\quad
v_2 = \begin{pmatrix}
-pd & 0 & 0 \\
0 & -pd & 0 \\
1 & 0 & 0 \\
\end{pmatrix},\quad
v_3 = \begin{pmatrix}
0 & -pd & 0 \\
0 & p & 0 \\
0 & 1 & 0 \\
\end{pmatrix}.
$$
If $d,p\neq0$, we apply $\psi$ with $\beta = \gamma = \varepsilon = 0$,
$\alpha pd = -1$, and $\delta d = 1$ and get
$$
\delta \psi(v_1) = \begin{pmatrix}
0 & 0 & 0 \\
1 & 1 & 0 \\
0 & 0 & 0 \\
\end{pmatrix},\quad
\alpha \psi(v_2) = \begin{pmatrix}
1 & 0 & 0 \\
0 & 1 & 0 \\
1 & 0 & 0 \\
\end{pmatrix},\quad
(\alpha/\delta)\psi(v_3) = \begin{pmatrix}
0 & 1 & 0 \\
0 & -1 & 0 \\
0 & 1 & 0 \\
\end{pmatrix},
$$
it is (T5).

If $d = 0$ and $p\neq0$, we apply $\psi$
with $\alpha p/\delta = 1$ and $\beta = \gamma = \varepsilon = 0$
and get $\psi(S) = \Span\{e_{21},e_{31},e_{22}+e_{32}\}$,
which is~(T2) up to the action of~$\Theta_{13}$ and transpose.

If $d\neq0$ and $p = 0$, we apply $\psi$ with
$\beta = \gamma = \varepsilon = 0$ and $\delta d = 1$ to get
$\psi(S) = \Span\{e_{21}+e_{22},e_{31},e_{32}\}$, it is (T2).

If $d = p = 0$, we have $\psi(S) = \Span\{e_{21},e_{31},e_{32}\}$, it is (T1).

{\sc Case IB}. $m = p\neq0$.
We have the second solution
$$
v_1 = \begin{pmatrix}
0 & 0 & 0 \\
1 & d & 0 \\
0 & 0 & 0 \\
\end{pmatrix},\quad
v_2 = \begin{pmatrix}
-pd & 0 & 0 \\
0 & -pd & 0 \\
1 & 0 & 0 \\
\end{pmatrix},\quad
v_3 = \begin{pmatrix}
p & 0 & 0 \\
0 & p & 0 \\
0 & 1 & 0 \\
\end{pmatrix},\quad p\neq0.
$$
If $d\neq0$, then we apply $\psi$ with $\beta = \gamma = \varepsilon = 0$,
$\alpha pd = -1$, and $\delta d = 1$ and get (T4).

If $d = 0$, then we apply $\psi$ with $\beta = \gamma = \varepsilon = 0$
and $\alpha p/\delta = 1$ to get
$\psi(S) = \Span\{e_{21},e_{31},e_{11}+e_{22}+e_{32}\}$,
which is~(T3) by the action of~$\Theta_{13}$ and transpose.

{\sc Case II}. $f = d\neq0$.

Then $m = 0$, and we have the third solution
$$
v_1 = \begin{pmatrix}
0 & 0 & -pd^2 \\
1 & d & pd \\
0 & 0 & d \\
\end{pmatrix},\quad
v_2 = \begin{pmatrix}
-pd & 0 & 0 \\
0 & -pd & 0 \\
1 & 0 & 0 \\
\end{pmatrix},\quad
v_3 = \begin{pmatrix}
0 & -dp & 0 \\
0 & p & 0 \\
0 & 1 & 0 \\
\end{pmatrix},\quad d\neq0.
$$
If $p\neq0$, then we apply $\psi$ with $\beta = \gamma = \varepsilon = 0$,
$\alpha pd = -1$, and $\delta d = 1$ and get (T6).

If $p = 0$, then we apply $\psi$ with $\beta = \gamma = \varepsilon = 0$
and $\delta d = 1$ to get (T3).

Theorem is proved.\hfill $\square$

{\bf Remark 3}.
All cases (T1)--(T6) from Theorem~3 lie in different orbits under action automorphisms
or antiautomorphisms of $M_3(\mathbb{C})$ preserving~$M$.
Indeed, the subalgebra~$S$ from (T1) is unique nilpotent one.
The cases (T2) and (T3) are the only cases which have one-dimensional semisimple part
in the Wedderburn---Malcev decomposition. They are not isomorphic, since
the subalgebra $S$ from (T2) but not from~(T3) has nonzero annihilator containing $e_{31}$.
There is a one-dimensional radical in~(T4)--(T6).
Note that the subalgebra~$S$ in the case~(T5) but not in~(T4) and (T6) is unital,
the matrix $e_{11}+e_{22}+e_{31}$ is its unit.
Further, the subalgebras~$S_4$ and $S_6$ from the cases (T4) and (T6) respectively are antiisomoprhic.
One can check directly that there are no $\psi\in\Aut(M_3(\mathbb{C}))$ preserving~$M$
such that neither $\varphi(S_4^T) = S_6$ nor $\psi(S_4) = S_6^T$.

\section{Unital (5,4)-decompositions}

Let us describe all 5-dimensional subalgebras in $M_3(\mathbb{C})$.

{\bf Lemma 5}.
Every 5-dimensional subalgebra in $M_3(\mathbb{C})$ up to action of $\Aut( M_3(\mathbb{C}) )$
and up to transpose is one of the following ones,

1) $\Span\{e_{11},e_{22},e_{23},e_{32},e_{33}\}$;

2) $\Span\{e_{11},e_{12},e_{13},e_{22},e_{33}\}$;

3) $\Span\{e_{11},e_{12},e_{13},e_{23},e_{22}+e_{33}\}$;

4) $\Span\{e_{11},e_{12},e_{13},e_{22},e_{23}\}$;

5) $\Span\{e_{11},e_{12},e_{13},e_{23},e_{33}\}$;

6) $\Span\{e_{11}+e_{33},e_{12},e_{13},e_{22},e_{23}\}$.

{\sc Proof}.
Let $S$ be a 5-dimensional subalgebra in $M_3(\mathbb{C})$.
By~\cite{AgoreMax,IovanovSistko}, we may assume that
$S \subset M$, where $M = \Span\{e_{11},e_{12},e_{13},e_{22},e_{23},e_{32},e_{33}\}$.
Denote the subalgebra $\Span\{e_{22}$, $e_{23},e_{32},e_{33}\}$ by~$N$.
We consider different cases of $\dim(S\cap N)$.

I) $\dim(S\cap N) = 4$.
It means that $N\subset S$.
It is easy to verify that the only possible case is 1).
Indeed, let $x = \alpha e_{11}+\beta e_{12}+\gamma e_{13}\in S$
with $(\beta,\gamma)\neq(0,0)$.
Up to $\Theta_{23}$, we may assume that $\beta\neq0$.
So, $xe_{22} = \beta e_{12}\in S$, i.e., $e_{12}$ lies in~$S$
as well $e_{12}e_{23} = e_{13}$, a~contradiction to the condition that $\dim S = 5$.

II) $\dim(S\cap N) = 3$.
Let us show that we may assume that $S\cap N = \Span\{e_{22},e_{23},e_{33}\}$.
Indeed, by~\cite{AgoreMax}, we have either $S\cap N = \Span\{e_{22},e_{23},e_{33}\}$
or $S\cap N = \Span\{e_{22},e_{32},e_{33}\}$.
Applying~$\Theta_{23}$ in the second case, we have stated the claim.
Let us clarify how we extend an automorphism $\varphi$ of~$N$
on the whole $M_3(\mathbb{C})$ preserving the property $\varphi(S)\subset M$.
Let $W\in M_2(\mathbb{C})$ be such that
$\varphi(X) = W^{-1}XW$ for every $X\in M_2(\mathbb{C})$,
then extend~$\varphi$ to an automorphism of $M_3(\mathbb{C})$ as a conjugation with the matrix
$\begin{pmatrix}
1 & 0 \\
0 & W \\
\end{pmatrix}$.

Analogously to the case~I, it is easy to see that we have exactly two possibilities:
either $S = \Span\{e_{12}$, $e_{13},e_{22},e_{23},e_{33}\}$
or $S = \Span\{e_{11},e_{13},e_{22},e_{23},e_{33}\}$.
Applying $\Theta_{13}\circ T$, we get~4) and~2) respectively.

III) $\dim(S\cap N) = 2$.

IIIa) $T = S\cap N$ is unital in~$N$.
Then we may assume that either
$T = \Span\{e_{22},e_{33}\}$ or
$T = \Span\{e_{22}+e_{33},e_{23}\}$.
In the first case it is 2) by the same reasons which are stated in~I.
In the second one, we have 3).

IIIb) $T = S\cap N$ is not unital in~$N$.
By~\cite{Meshulam} and up to the action of $\Theta_{23}$,
we have two cases:
$T = \Span\{e_{22},e_{23}\}$ or
$T = \Span\{e_{23},e_{33}\}$.
Let $T = \Span\{e_{22},e_{23}\}$.
We want to state that we may assume that $e_{12},e_{13}\in S$.
Let $x = \begin{pmatrix}
a & b & c \\
0 & 0 & 0 \\
0 & d & e
\end{pmatrix}\in S$ with $b\neq0$.
Then $xe_{22} = be_{12} + de_{32}\in S$ and
$xe_{23} = be_{13} + de_{33}\in S$.
Suppose that $d\neq0$.
Applying the automorphism~$\varphi$~\eqref{AutPreserv} which preserves~$M$
with $\beta = \lambda = \mu = 0$ and $\nu b = \gamma d$, we get
$\varphi(S)\ni e_{22},e_{23},e_{32},e_{33}$,
it is case I.

So, $S = \Span\{e_{12},e_{13},e_{22},e_{23},\alpha e_{11}+\beta e_{33}\}$.
Since $S$ is 5-dimensional, we have either $\beta = 0$ or $\alpha = \beta$ (the case
$\alpha = 0$ is impossible, otherwise $e_{33}\in S\cap N$).
In the first case, it is 4), in the second one it is 6).

When $T = \Span\{e_{23},e_{33}\}$, we have $e_{13}\in S$.
Note that $S$ has no nonzero projection on $e_{32}$, otherwise
$e_{33}(e_{32}+\alpha e_{33}+\ldots) = e_{32} + \alpha e_{33}\in S$,
i.\,e. $e_{32}\in S$, a contradiction.
So, $V = S\cap \Span\{e_{11},e_{12},e_{22}\}$
is a~2-dimensional subalgebra in $K = \Span\{e_{11},e_{12},e_{21},e_{22}\}$.
As above, we may extend an automorphism $\varphi$ of~$K$
on the whole $M_3(\mathbb{C})$ preserving the property that $e_{13},e_{23},e_{33}\in\varphi(S)$.

If $V$~is not unital in~$K$, then up to action of $\Theta_{12}$ we have either
$V = \Span\{e_{11},e_{12}\}$ (it is 5)) or $V = \Span\{e_{12},e_{22}\}$
(applying $\Theta_{13}\circ T$, we get 4)).

If $V$~is unital in~$K$, then up to action of $\Theta_{12}$, we have either
$V = \Span\{e_{11},e_{22}\}$ or $V = \Span\{e_{11}+e_{22},e_{12}\}$.
Applying $\Theta_{13}\circ T$, we get 2) and 3) respectively.
\hfill $\square$

{\bf Remark 4}.
All subalgebras listed in Lemma~5 are pairwise neither isomorphic nor antiisomorphic.
Indeed, the subalgebra 1) is unique semisimple one from the list.
Further, the subalgebra from 2) is unique which has two-dimensional radical.
The subalgebras 4) and 5) are unique non-unital one, moreover,
there exists an idempotent in the subalgebra from 4) but not from 5) such that
its products with the second power of the radical is zero.
The last argument distinguishes subalgebras from 3) and 6) as well.

\subsection{Case of unital 4-dimensional subalgebra}

Let $M$ be a non-unital 5-dimensional algebra,
by Lemma~5 we may assume up to transpose and action of $\Aut(M_3(\mathbb{C}))$
that either $M = \Span\{e_{11},e_{12},e_{13},e_{22},e_{23}\}$ or
$M = \Span\{e_{11},e_{12},e_{13},e_{23},e_{33}\}$.
Then the group $\Aut(M_3(\mathbb{C}))$ preserving~$M$
coincides with the group of automorphisms of the subalgebra of upper-triangular matrices.

{\bf Theorem 4}.
Every direct decomposition of $M_3(\mathbb{C})$
with two subalgebras of the dimensions~4 and~5,
where a 4-dimensional subalgebra is unital,
up to transpose and up to action of $\Aut(M_3(\mathbb{C}))$
is isomorphic to $S\oplus M$, where either
$M = \Span\{e_{11},e_{12},e_{13},e_{22},e_{23}\}$ or
$M = \Span\{e_{11},e_{12},e_{13},e_{23},e_{33}\}$, and
where $S$~is one of the following subalgebras:

(U1) $S = \Span\{E,e_{21},e_{31},e_{32}\}$;

(U2) $S = \Span\{E,e_{21}+e_{22},e_{31},e_{32}\}$;

(U3) $S = \Span\{E,e_{21},e_{31},e_{32}+e_{22}\}$;

(U4) $S = \Span\{E,e_{21}+e_{22},e_{31}+e_{13}-e_{23},e_{11}+e_{12}+e_{31}+e_{32}\}$;

(U5) $S = \Span\{E,e_{21},e_{31},e_{11}+pe_{22}+(1-p)e_{23}+e_{32}\}$;

(U6) $S = \Span\{E,e_{21}-e_{23},e_{31}+e_{11}-e_{13}-e_{33},e_{11}+e_{12}-e_{13}+pe_{22}+(1-p)e_{23}+e_{32}\}$;

(U7) $S = \Span\{E,e_{12}+e_{22},e_{31}+m(m-1)(e_{13}-e_{23})-e_{33},
e_{31}+e_{32}+(m+1)(e_{11}+e_{12})\}$;

(U8) $S = \Span\{E,
e_{12}+e_{22}-e_{32},
(m-p+1)(e_{11}+e_{12}-e_{13})
+e_{31}+e_{32}-e_{33},
e_{11}+e_{31}+(m-1)(m-p+1)e_{13}-m(m-p)e_{23}+(p-1)e_{33}\}$;

{\sc Proof}.
The subalgebra~$S$ is generated by $E$ and by matrices
$$
v_1 = \begin{pmatrix}
a & b & c \\
1 & d & e \\
0 & 0 & 0 \\
\end{pmatrix},\quad
v_2 = \begin{pmatrix}
g & h & i \\
0 & j & k \\
1 & 0 & 0 \\
\end{pmatrix},\quad
v_3 = \begin{pmatrix}
m & n & s \\
0 & p & q \\
0 & 1 & 0 \\
\end{pmatrix}.
$$
Applying $\psi$~\eqref{AutU}
with $\gamma = \varepsilon = 0$ and $\beta = a\delta$,
we may assume that $\underline{a = 0}$.

From
$$
v_1^2 - dv_1 - bE = \begin{pmatrix}
0 & 0 & be-cd \\
0 & 0 & c \\
0 & 0 & -b
\end{pmatrix},
$$
it follows that $b = c = 0$.

From
$$
v_3^2 - pv_3 - qE = \begin{pmatrix}
m^2-mp-q & mn+s & nq+ms-ps \\
0 & 0 & 0 \\
0 & 0 & 0
\end{pmatrix},
$$
we conclude that $s = -mn$ and $q = m(m-p)$.

From the equality
$$
v_2^2 -gv_2 - kv_1 - hv_3 - iE = \begin{pmatrix}
-hm & h(j-n) & h(k+mn) \\
0 & j^2 - jg -dk - hp - i & k(j-e-g) - hq \\
0 & 0 & 0
\end{pmatrix},
$$
we derive
$$
hm = hk = 0, \quad
h(j-n) = 0,\quad
j^2 - jg -dk - hp - i = 0, \quad
k(j-e-g) - hq = 0.
$$

From
$$
v_3v_1 - v_2 - pv_1 - dv_3 = \begin{pmatrix}
n-g-dm & -h & dmn+ne-i \\
0 & -j-dp & dmp-k-dm^2 \\
0 & 0 & e \\
\end{pmatrix},
$$
we obtain $h = 0$, and the system of equations
$$
n = e + g + dm, \quad
j = - e - dp, \quad
k = dm(p-m),\quad
i = n(n-g).
$$

From the product $v_1v_3$, we get the equality
$n+e+d(p-m) = 0$.

From
$$
v_1v_2 - (g+e)v_1 = \begin{pmatrix}
0 & 0 & 0 \\
0 & d(j-e-g) & dk - e^2 - ge + i \\
0 & 0 & 0 \\
\end{pmatrix},
$$
we get the equations
$$
d(j-e-g) = 0, \quad
dk - e^2 - ge + i = 0.
$$
Finally, we have
\begin{gather*}
v_3 v_2 - jv_3 - m(m-p)v_1 = \begin{pmatrix}
m(g-j-n) & 0 & im + kn + jmn \\
0 & dm(p-m) &  m(j + e)(p - m) + kp \\
0 & 0 & k
\end{pmatrix}, \\
v_2 v_3 - mv_2 - nv_3 + mn E = \begin{pmatrix}
0 & gn-n^2+i & m(n^2-gn-i) \\
0 & k+(n-j)(m-p) & m((n-j)(p-m)-k) \\
0 & 0 & 0
\end{pmatrix}.
\end{gather*}

Summing up, we have
\begin{gather*}
n = d(m-p)-e, \quad
j = - e - dp, \quad
k = dm(p-m),\\
g = -2e-dp, \quad
i = (d(m-p)-e)(dm+e),
\end{gather*}
and
\begin{gather*}
v_1 = \begin{pmatrix}
0 & 0 & 0 \\
1 & d & e \\
0 & 0 & 0 \\
\end{pmatrix},\quad
v_2 = \begin{pmatrix}
-2e-dp & 0 & (dm+e)(d(m-p)-e) \\
0 & -e-dp & -dm(m-p) \\
1 & 0 & 0 \\
\end{pmatrix}, \\
v_3 = \begin{pmatrix}
m & d(m-p)-e & -m(d(m-p)-e) \\
0 & p & m(m-p) \\
0 & 1 & 0 \\
\end{pmatrix}.
\end{gather*}

{\sc Case I}: $e\neq0$.

{\sc Case IA}: $d\neq0$.
We apply $\psi$ with $\beta = \gamma = \varepsilon = 0$, $d\delta = -\alpha e = 1$
and get
\begin{gather*}
v_1 = \begin{pmatrix}
0 & 0 & 0 \\
1 & 1 & -1 \\
0 & 0 & 0 \\
\end{pmatrix},\quad
v_2 = \begin{pmatrix}
2-p & 0 & (m-1)(m-p+1) \\
0 & 1-p & -m(m-p) \\
1 & 0 & 0 \\
\end{pmatrix}, \\
v_3 = \begin{pmatrix}
m & m-p+1 & -m(m-p+1) \\
0 & p & m(m-p) \\
0 & 1 & 0 \\
\end{pmatrix}.
\end{gather*}
It is~(U8).

\newpage
{\sc Case IB}: $d = 0$.

{\sc Case IBA}: $m = 0$.
If $p = 0$, we apply $\psi$ with $\beta = \gamma = \varepsilon = 0$ and $\alpha e = -1$
and get the subalgebra~$\psi(S) = \Span\{E,e_{21}-e_{23},e_{31}+e_{11}-e_{13}-e_{33},e_{32}+e_{12}\}$.
Acting by~$\psi$ with $\beta = \varepsilon = 0$, $\gamma = \alpha$, we get~(U1).
If $p\neq0$, we use $\psi$ with $\beta = \gamma = \varepsilon = 0$,
$\alpha e = -1$, and $\delta = \alpha p$, so, we have the subalgebra
$\psi(S) = \Span\{E,e_{21}-e_{23},e_{31}+e_{11}-e_{13}-e_{33},e_{32}+e_{12}+e_{22}\}$.
After action of~$\psi$ defined with $\beta = \varepsilon = 0$, $\gamma = \delta = \alpha$, we get~(U3).

{\sc Case IBB}: $m\neq0$.
We apply $\psi$ with $\beta = \gamma = \varepsilon = 0$,
$\alpha e = -1$, and $\delta = \alpha m$, we get
$$
v_1 = \begin{pmatrix}
0 & 0 & 0 \\
1 & 0 & -1 \\
0 & 0 & 0 \\
\end{pmatrix},\quad
v_2 = \begin{pmatrix}
2 & 0 & -1 \\
0 & 1 & 0 \\
1 & 0 & 0 \\
\end{pmatrix}, \quad
v_3 = \begin{pmatrix}
1 & 1 & -1 \\
0 & p & 1-p \\
0 & 1 & 0 \\
\end{pmatrix},
$$
it is~(U6).

{\sc Case II}: $e = 0$.

{\sc Case IIA}: $d = 0$.

{\sc Case IIAA}: $m = 0$.
If $p = 0$, then we have~(U1).
If $p\neq0$, then
we apply $\psi$ with $\beta = \gamma = \varepsilon = 0$ and $\delta = \alpha p$,
thus $\psi(S) = \Span\{E,e_{21},e_{31},e_{22}+e_{32}\}$, it is~(U3).

{\sc Case IIAB}: $m\neq0$.
We apply $\psi$ with $\beta = \gamma = \varepsilon = 0$ and $\delta = \alpha m$
and get~(U5).

{\sc Case IIB}: $d\neq0$.

{\sc Case IIBA}: $p = 0$.
If $m = 0$, then we apply $\psi$
with $\beta = \gamma = \varepsilon = 0$ and $\delta d = 1$ and get~(U2).
If $m\neq0$, then we apply $\psi$
with $\beta = \gamma = \varepsilon = 0$, $\delta d = 1$,
and $\alpha dm = 1$. So, we obtain~(U4).

{\sc Case IIBB}: $p\neq0$.
We apply $\psi$ with $\beta = \gamma = \varepsilon = 0$, $\delta d = 1$,
and $\alpha dp = -1$. Therefore, we may assume that $S$ is generated by $E$ and by
$$
v_1 = \begin{pmatrix}
0 & 0 & 0 \\
1 & 1 & 0 \\
0 & 0 & 0 \\
\end{pmatrix}, \quad
v_2 = \begin{pmatrix}
1 & 0 & m(m+1) \\
0 & 1 & -m(m+1) \\
1 & 0 & 0 \\
\end{pmatrix}, \quad
v_3 = \begin{pmatrix}
m & m+1 & -m(m+1) \\
0 & -1 & m(m+1) \\
0 & 1 & 0 \\
\end{pmatrix},
$$
it is~(U7).

Theorem is proved.\hfill $\square$

{\bf Remark 5}.
Cases~(U2) and~(U3) coincide by the action of~$\Theta_{13}\circ T$,
when $M = \Span\{e_{11},e_{12},e_{13},e_{23},e_{33}\}$.

\subsection{Case of unital 5-dimensional subalgebra. I}

Let $M = \Span\{e_{11},e_{22},e_{23},e_{32},e_{33}\}$.
It is easy to show that an automorphism $\chi\in\Aut(M_3(\mathbb{C}))$
preserving~$M$ has the form~\eqref{AutPreserv} with $\beta = \gamma = 0$.

{\bf Theorem 5}.
Consider a direct decomposition of $M_3(\mathbb{C})$
with two subalgebras of the dimensions~4 and~5,
where a 5-dimensional subalgebra up to transpose and action of $\Aut(M_3(\mathbb{C}))$
is isomorphic to $M = \Span\{e_{11},e_{22},e_{23},e_{32},e_{33}\}$.
Then such decomposition
up to transpose and up to action of $\Aut(M_3(\mathbb{C}))$
is isomorphic to $S\oplus M$,
where $S$~is one of the following subalgebras:

(V1) $S = \Span\{e_{21}-e_{23},e_{11}+e_{22}+e_{31},e_{12}+e_{32},e_{13}+e_{22}+e_{33}\}$;

(V2) $S = \Span\{e_{21}-e_{23},e_{11}+e_{22}+e_{31},e_{12}+e_{23}+e_{32},e_{13}+e_{22}+e_{33}\}$;

(V3) $S = \Span\{e_{21}-e_{23},e_{11}+e_{31},e_{12}+e_{32},e_{13}+e_{33}\}$;

(V4) $S = \Span\{e_{21}+e_{22}-e_{23},e_{11}+e_{31},e_{12}+e_{32},e_{13}+e_{33}\}$;

(V5) $S = \Span\{e_{21}-e_{23},e_{11}+e_{22}+e_{31},e_{12}+e_{22}+se_{23}+e_{32},e_{13}+e_{22}+e_{33}\}$;

(V6) $S = \Span\{e_{21}+be_{22}-(b+1)e_{23},e_{11}+e_{21}+e_{31},e_{12}+e_{22}+e_{32},e_{13}+e_{23}+e_{33}\}$.

{\sc Proof}.
The subalgebra~$S$ is generated by matrices
$$
v_1 = \begin{pmatrix}
a & 0 & 0 \\
1 & b & c \\
0 & d & e \\
\end{pmatrix},\quad
v_2 = \begin{pmatrix}
g & 0 & 0 \\
0 & h & i \\
1 & j & k \\
\end{pmatrix},\quad
v_3 = \begin{pmatrix}
m & 1 & 0 \\
0 & n & s \\
0 & p & q \\
\end{pmatrix},\quad
v_4 = \begin{pmatrix}
u & 0 & 1 \\
0 & v & x \\
0 & y & z \\
\end{pmatrix}.
$$

Analyzing the $e_{11}$-coordinate of the products $v_i v_j$,
where $i=1,2$, $j=3,4$, we get
$$
am = au = gm = gu = 0.
$$
Up to the action of $\Theta_{23}$ and up to transpose,
we may assume that $a = 0$.
If $g\neq 0$, then $u = m = 0$.
If $g = 0$, then applying $\varphi$~\eqref{AutPreserv}
with $\beta = \gamma = 0$ and corresponding $\kappa,\lambda,\mu,\nu$,
we may assume that one of $m,u$ is zero.
If only one of $m,u$ is zero, this case
is equivalent to the case $a=m=u=0$ and $g\neq0$.
The case when all $a,g,m,u$ are zero is impossible
by the analysis of the $e_{11}$-coordinate
$1+am-an-pg-bm-cu$ from the equality $v_3v_1 = nv_1 + pv_2 + bv_3 + cv_4$.

So, we assume that $\underline{a = m = u = 0}$ and $\underline{g\neq0}$.
From
$$
v_1^2 = bv_1 + dv_2, \quad
v_1v_2 = (g+c)v_1 + ev_2,\quad
v_2v_1 = hv_1 + jv_2,\quad
v_2^2 = iv_1 + (g+k)v_2,
$$
we get $dg = eg = jg = kg = 0$, it means that $d = e = j = k = 0$.

The $e_{11}$-coordinates of the products $v_3v_1$  and $v_2v_4$
give $p = z = 1/g$. So, we have
$$
v_1 = \begin{pmatrix}
0 & 0 & 0 \\
1 & b & c \\
0 & 0 & 0 \\
\end{pmatrix},\quad
v_2 = \begin{pmatrix}
g & 0 & 0 \\
0 & h & i \\
1 & 0 & 0 \\
\end{pmatrix},\quad
v_3 = \begin{pmatrix}
0 & 1 & 0 \\
0 & n & s \\
0 & 1/g & q \\
\end{pmatrix},\quad
v_4 = \begin{pmatrix}
0 & 0 & 1 \\
0 & v & x \\
0 & y & 1/g \\
\end{pmatrix}.
$$
From the equalities
\begin{gather*}
v_1v_3 = 0, \quad
v_3v_1 = nv_1 + (1/g)v_2 + bv_3 + cv_4,\quad
v_1v_4 = 0,\quad
v_4v_1 = vv_1 + yv_2,\\
v_2v_3 = gv_3,\quad
v_3v_2 = sv_1 + qv_2 + hv_3 + iv_4,\quad
v_2v_4 = gv_4,\quad
v_4v_2 = xv_1 + (1/g)v_2,
\end{gather*}
we get $y = q = 0$ and obtain the following system of equations
\begin{equation}\label{5-4-I:MainSystem}
\begin{gathered}
\{b,h-g\}\times\{v,s,n-x\} = 0, \\
1+bn+c/g = 0, \quad
bn+h/g+cv = 0, \\
cx + i/g = 0, \quad
x(c+g-h) = 0, \quad
iv = 0, \\
i(n-x) - s(c+h) = 0, \quad
v(c+h) = 0.
\end{gathered}
\end{equation}

{\sc Case I}: $b = 0$ and $h = g\neq0$.
Then~\eqref{5-4-I:MainSystem} gives
$c = -g$, $h = g^2v$, it means that $v = 1/g$.
Thus, $i = x = 0$. It is the first solution with $g\neq0$
$$
v_1 = \begin{pmatrix}
0 & 0 & 0 \\
1 & 0 & -g \\
0 & 0 & 0 \\
\end{pmatrix},\quad
v_2 = \begin{pmatrix}
g & 0 & 0 \\
0 & g & 0 \\
1 & 0 & 0 \\
\end{pmatrix},\quad
v_3 = \begin{pmatrix}
0 & 1 & 0 \\
0 & n & s \\
0 & 1/g & 0 \\
\end{pmatrix},\quad
v_4 = \begin{pmatrix}
0 & 0 & 1 \\
0 & 1/g & 0 \\
0 & 0 & 1/g \\
\end{pmatrix}.
$$
If $n = s = 0$, then we apply $\varphi$~\eqref{AutPreserv} with
$\beta = \gamma = \lambda = \mu = 0$ and $\nu g = 1$ to get~(V1).
If $n = 0$ and $s\neq0$, then we apply $\varphi$ with
$\beta = \gamma = \lambda = \mu = 0$, $\nu g = 1$, and $\kappa^2 = s/g$ to get~(V2).
If $n\neq0$, then we apply $\varphi$ with
$\beta = \gamma = \lambda = \mu = 0$, $\nu g = 1$, and $\kappa = n$ to get~(V5).

{\sc Case II}: $(b,h-g) \neq (0,0)$.
Then $v = s = 0$ and $n = x$.
The system~\eqref{5-4-I:MainSystem} is equivalent to the following one,
$$
1 + bn + c/g = 0, \quad
bn + h/g = 0, \quad
cn + i/g = 0, \quad
n(c+g-h) = 0.
$$
Thus, $h = -bng$, $c = h-g = -g(bn+1)$, $i = -cng = ng^2(bn+1)$.
It is the second solution with $g\neq0$
\begin{gather*}
v_1 = \begin{pmatrix}
0 & 0 & 0 \\
1 & b & -g(bn+1) \\
0 & 0 & 0 \\
\end{pmatrix},\quad
v_2 = \begin{pmatrix}
g & 0 & 0 \\
0 & -bng & ng^2(bn+1) \\
1 & 0 & 0 \\
\end{pmatrix}, \\
v_3 = \begin{pmatrix}
0 & 1 & 0 \\
0 & n & 0 \\
0 & 1/g & 0 \\
\end{pmatrix},\quad
v_4 = \begin{pmatrix}
0 & 0 & 1 \\
0 & 0 & n \\
0 & 0 & 1/g \\
\end{pmatrix}.
\end{gather*}
If $b = n = 0$, then we apply $\varphi$ with
$\beta = \gamma = \lambda = \mu = 0$ and $\nu g = 1$ to get~(V3).
If $n = 0$ and $b\neq0$, then we apply $\varphi$ with
$\beta = \gamma = \lambda = \mu = 0$, $\nu g = 1$, and $\kappa b = 1$ to get~(V4).
If $n\neq0$, then we apply $\varphi$ with
$\beta = \gamma = \lambda = \mu = 0$, $\nu g = 1$, and $\kappa = n$ to get~(V6).

Theorem is proved.\hfill $\square$

\subsection{Case of unital 5-dimensional subalgebra. II}

Let $M = \Span\{e_{11},e_{12},e_{13},e_{23},e_{22}+e_{33}\}$.
Note that all automorphisms of the upper-triangular matrices preserve $M$.
It is not difficult to show the inverse, it means that
an automorphism $\chi\in\Aut(M_3(\mathbb{C}))$ preserving~$M$
is defined by~\eqref{AutU}. Indeed, since an
automorphism~$\chi\in\Aut(M_3(\mathbb{C}))$ has to preserve the radical of~$M$,
we have $\chi(e_{13}) = \alpha e_{13}$, $\alpha\neq0$,
$\chi(e_{12}) = \beta e_{12} + \gamma e_{13} + \delta e_{23}$.
Since $e_{12}^2 = 0$, we have $\beta\delta = 0$. If $\beta = 0$, then
$\chi(e_{23}) = \epsilon e_{12} + \tau e_{13}$ with $\epsilon\neq0$.
Thus, $\chi(e_{13}) \neq \chi(e_{12})\chi(e_{23})$, a contradiction.
Hence, $\delta = 0$, and $\chi$ satisfies~\eqref{AutPreserv} with $\mu = 0$,
it is~\eqref{AutU}.

{\bf Theorem 6}.
Consider a direct decomposition of $M_3(\mathbb{C})$
with two subalgebras of the dimensions~4 and~5,
where a 5-dimensional subalgebra up to transpose and action of $\Aut(M_3(\mathbb{C}))$
is isomorphic to $M = \Span\{e_{11},e_{12},e_{13},e_{23},e_{22}+e_{33}\}$.
Then such decomposition up to transpose and up to action
of $\Aut(M_3(\mathbb{C}))$ is isomorphic to $S\oplus M$,
where $S$~is one of the following subalgebras:

(X1) $S = \Span\{e_{21},e_{22},e_{31},e_{32}\}$;

(X2) $S = \Span\{e_{21},e_{31},e_{32},e_{33}\}$;

(X3) $S = \Span\{e_{21},e_{31},e_{32},e_{11}+e_{22}\}$;

(X4) $S = \Span\{e_{21},e_{31},e_{32},e_{11}+e_{33}\}$;

(X5) $S = \Span\{e_{11}+e_{31},e_{12}+e_{32},e_{21},e_{22}\}$;

(X6) $S = \Span\{e_{21}+e_{12},e_{31},e_{32},e_{11}+e_{22}\}$;

(X7) $S = \Span\{e_{21}+e_{22},e_{31},e_{32},e_{33}\}$.

{\sc Proof}. The subalgebra~$S$ is generated by matrices
$$
v_1 = \begin{pmatrix}
a & b & c \\
1 & d & e \\
0 & 0 & d \\
\end{pmatrix}, \quad
v_2 = \begin{pmatrix}
g & h & i \\
0 & j & k \\
1 & 0 & j \\
\end{pmatrix}, \quad
v_3 = \begin{pmatrix}
m & n & s \\
0 & p & q \\
0 & 1 & p \\
\end{pmatrix}, \quad
v_4 = \begin{pmatrix}
u & v & x \\
0 & y & z \\
0 & 0 & y+1 \\
\end{pmatrix}.
$$
Applying~$\psi$~\eqref{AutU} with $\beta = 0$,
$\gamma = -\alpha e$, and $\varepsilon = \alpha z$,
we may assume that $\underline{e = z = 0}$.

From the products
\begin{gather*}
v_1^2 = (a+d)v_1 - bv_4, \quad
v_1v_2 = gv_1 + dv_2 - hv_4, \quad
v_1v_3 = mv_1 + dv_3 - nv_4, \\
v_1v_4 = uv_1 + (d-v)v_4, \
v_2v_1 = jv_1 + av_2 + bv_3 + cv_4, \
v_2^2 = kv_1 + (g+j)v_2+ hv_3+iv_4, \\
v_2v_3 = mv_2 + (j+n)v_3+(s-k)v_4, \quad
v_2v_4 = uv_2 + vv_3 + (j+x)v_4, \\
v_3v_1 = pv_1 + v_2 + dv_3, \quad
v_3v_2 = qv_1 + pv_2 + jv_3 + kv_4, \quad
v_3^2  = 2pv_3, \quad
v_3v_4 = yv_3 + pv_4, \\
v_4v_1 = yv_1 + dv_4, \quad
v_4v_2 = (y+1)v_2 + jv_4, \quad
v_4v_3 = (y+1)v_3 + pv_4, \quad
v_4^2 = (2y+1)v_4,
\end{gather*}
we get the system of 80 equations, which radical computed by Singular~\cite{Singular} gives
$c = h = i = j = k = m = p = q = s = x = 0$, $n = g$ and the system of 11~equations
\begin{gather*}
bg = 0, \quad
dg = 0, \quad
dy = 0, \quad
gu = 0, \allowdisplaybreaks \\
b(u-y) = 0, \quad
y(y+1) = 0, \quad
g(y+1) = 0, \quad
u(u-2y-1) = 0, \\
du-v(y+1) = 0, \quad
a(u-y)-du+v = 0, \quad
ad - b(u+1) = 0.
\end{gather*}

{\sc Case I}: $g\neq0$.
Then $y = -1$, $b = d = u = 0$, and $v = -a$. It is the first solution
$$
v_1 = \begin{pmatrix}
a & 0 & 0 \\
1 & 0 & 0 \\
0 & 0 & 0 \\
\end{pmatrix}, \quad
v_2 = \begin{pmatrix}
g & 0 & 0 \\
0 & 0 & 0 \\
1 & 0 & 0 \\
\end{pmatrix}, \quad
v_3 = \begin{pmatrix}
0 & g & 0 \\
0 & 0 & 0 \\
0 & 1 & 0 \\
\end{pmatrix}, \quad
-v_4 = \begin{pmatrix}
0 & a & 0 \\
0 & 1 & 0 \\
0 & 0 & 0 \\
\end{pmatrix},\ g\neq0.
$$

{\sc Case II}: $g = 0$.
We have the remaining system of the equations
\begin{gather*}
dy = 0, \quad
b(u-y) = 0, \quad
y(y+1) = 0, \quad
u(u-2y-1) = 0, \\
du-v(y+1) = 0, \quad
a(u-y)-du+v = 0, \quad
ad - b(u+1) = 0.
\end{gather*}

{\sc Case IIA}: $y = -1$.
Then $d = 0$ and we have the system
$$
a(u+1)+v = 0, \quad
b(u+1) = 0, \quad
u(u+1) = 0.
$$

{\sc Case IIAA}: $u = -1$. Thus, $v = 0$,
it is the second solution
$$
v_1 = \begin{pmatrix}
a & b & 0 \\
1 & 0 & 0 \\
0 & 0 & 0 \\
\end{pmatrix}, \quad
v_2 = \begin{pmatrix}
0 & 0 & 0 \\
0 & 0 & 0 \\
1 & 0 & 0 \\
\end{pmatrix}, \quad
v_3 = \begin{pmatrix}
0 & 0 & 0 \\
0 & 0 & 0 \\
0 & 1 & 0 \\
\end{pmatrix}, \quad
-v_4 = \begin{pmatrix}
1 & 0 & 0 \\
0 & 1 & 0 \\
0 & 0 & 0 \\
\end{pmatrix}.
$$

{\sc Case IIAB}: $u = 0$. So, $b = 0$ and $v = -a$,
it is the third solution
$$
v_1 = \begin{pmatrix}
a & 0 & 0 \\
1 & 0 & 0 \\
0 & 0 & 0 \\
\end{pmatrix}, \quad
v_2 = \begin{pmatrix}
0 & 0 & 0 \\
0 & 0 & 0 \\
1 & 0 & 0 \\
\end{pmatrix}, \quad
v_3 = \begin{pmatrix}
0 & 0 & 0 \\
0 & 0 & 0 \\
0 & 1 & 0 \\
\end{pmatrix}, \quad
-v_4 = \begin{pmatrix}
0 & a & 0 \\
0 & 1 & 0 \\
0 & 0 & 0 \\
\end{pmatrix}.
$$

{\sc Case IIB}: $y = 0$. Then we have the system
$$
au = 0, \quad bu = 0, \quad u(u-1) = 0, \quad
v = du, \quad ad - b(u+1) = 0.
$$

{\sc Case IIBA}: $u = 0$. So, $v = 0$ and $b = ad$.
It is the fourth soultion
$$
v_1 = \begin{pmatrix}
a & ad & 0 \\
1 & d & 0 \\
0 & 0 & d \\
\end{pmatrix}, \quad
v_2 = \begin{pmatrix}
0 & 0 & 0 \\
0 & 0 & 0 \\
1 & 0 & 0 \\
\end{pmatrix}, \quad
v_3 = \begin{pmatrix}
0 & 0 & 0 \\
0 & 0 & 0 \\
0 & 1 & 0 \\
\end{pmatrix}, \quad
v_4 = \begin{pmatrix}
0 & 0 & 0 \\
0 & 0 & 0 \\
0 & 0 & 1 \\
\end{pmatrix}.
$$

{\sc Case IIBB}: $u = 1$. Then $a = b = 0$ and $v = d$.
It is the fifth soultion
$$
v_1 = \begin{pmatrix}
0 & 0 & 0 \\
1 & d & 0 \\
0 & 0 & d \\
\end{pmatrix}, \quad
v_2 = \begin{pmatrix}
0 & 0 & 0 \\
0 & 0 & 0 \\
1 & 0 & 0 \\
\end{pmatrix}, \quad
v_3 = \begin{pmatrix}
0 & 0 & 0 \\
0 & 0 & 0 \\
0 & 1 & 0 \\
\end{pmatrix}, \quad
v_4 = \begin{pmatrix}
1 & d & 0 \\
0 & 0 & 0 \\
0 & 0 & 1 \\
\end{pmatrix}.
$$

Let us gather obtained cases as follows,
where the first and the third cases come together in 1) without restrictions on~$g$,
\begin{align*}
& 1)\ v_1 = \begin{pmatrix}
a & 0 & 0 \\
1 & 0 & 0 \\
0 & 0 & 0 \\
\end{pmatrix},
&& v_2 = \begin{pmatrix}
g & 0 & 0 \\
0 & 0 & 0 \\
1 & 0 & 0 \\
\end{pmatrix},
&& v_3 = \begin{pmatrix}
0 & g & 0 \\
0 & 0 & 0 \\
0 & 1 & 0 \\
\end{pmatrix},
&& v_4 = \begin{pmatrix}
0 & a & 0 \\
0 & 1 & 0 \\
0 & 0 & 0 \\
\end{pmatrix}, \allowdisplaybreaks \\
& 2)\ v_1 = \begin{pmatrix}
a & b & 0 \\
1 & 0 & 0 \\
0 & 0 & 0 \\
\end{pmatrix},
&& v_2 = \begin{pmatrix}
0 & 0 & 0 \\
0 & 0 & 0 \\
1 & 0 & 0 \\
\end{pmatrix},
&& v_3 = \begin{pmatrix}
0 & 0 & 0 \\
0 & 0 & 0 \\
0 & 1 & 0 \\
\end{pmatrix},
&& v_4 = \begin{pmatrix}
1 & 0 & 0 \\
0 & 1 & 0 \\
0 & 0 & 0 \\
\end{pmatrix}, \\
& 3)\ v_1 = \begin{pmatrix}
a & ad & 0 \\
1 & d & 0 \\
0 & 0 & d \\
\end{pmatrix},
&& v_2 = \begin{pmatrix}
0 & 0 & 0 \\
0 & 0 & 0 \\
1 & 0 & 0 \\
\end{pmatrix},
&& v_3 = \begin{pmatrix}
0 & 0 & 0 \\
0 & 0 & 0 \\
0 & 1 & 0 \\
\end{pmatrix},
&& v_4 = \begin{pmatrix}
0 & 0 & 0 \\
0 & 0 & 0 \\
0 & 0 & 1 \\
\end{pmatrix}, \\
& 4)\ v_1 = \begin{pmatrix}
0 & 0 & 0 \\
1 & d & 0 \\
0 & 0 & d \\
\end{pmatrix},
&& v_2 = \begin{pmatrix}
0 & 0 & 0 \\
0 & 0 & 0 \\
1 & 0 & 0 \\
\end{pmatrix},
&& v_3 = \begin{pmatrix}
0 & 0 & 0 \\
0 & 0 & 0 \\
0 & 1 & 0 \\
\end{pmatrix},
&& v_4 = \begin{pmatrix}
1 & d & 0 \\
0 & 0 & 0 \\
0 & 0 & 1 \\
\end{pmatrix}.
\end{align*}

{\sc Case 1}.
Let us apply~$\psi$~\eqref{AutU} defined with
$\gamma = \varepsilon = 0$, $\alpha = 1$, and $\beta = \delta a$
to get $\psi(S) = \Span\{e_{21},e_{22},ge_{11}+e_{31},ge_{12}+e_{32}\}$.
If $g = 0$, then we have (X1).
If $g\neq0$, then we apply~$\psi$
with $\beta = \gamma = \varepsilon = 0$ and $\alpha g = 1$
to get~(X5).

{\sc Case 2}.
We apply~$\psi$ with $\gamma = \varepsilon = 0$ and $\beta = \delta a/2$ to get
$\psi(S) = \Span\{b'e_{12}+e_{21},e_{31},e_{32},e_{11}+e_{22}\}$, where $b' = \delta^2(b+a^2/4)$.
If $b' = 0$, we have~(X3).
If $b'\neq 0$, we apply ~$\psi$ with $\beta = \gamma = \varepsilon = 0$
and $\delta^2 b' = 1$ to get~(X6).

{\sc Case 3}.
We apply $\psi$ with $\gamma = \varepsilon = 0$, $\beta = \delta a$
to get $\psi(S) = \Span\{e_{31},e_{32},e_{33},e_{21}+\delta(a+d)e_{22}\}$.
If $a+d = 0$, then we have~(X2).
Otherwise, we get~(X7) when $\delta(a+d) = 1$.

{\sc Case 4}.
We apply $\psi$ with
$\gamma = \varepsilon = 0$ and $\beta = - \delta d$ to get~(X4).

Theorem is proved.\hfill $\square$

{\bf Remark 6}.
All cases (X1)--(X7) from Theorem~6 lie in different orbits under action of automorphisms
or antiautomorphisms of $M_3(\mathbb{C})$ preserving~$M$.
Indeed, case~(X5) is only one where the subalgebra is semisimple, it is isomorphic to~$M_2(\mathbb{C})$.
The radical of the subalgebra~$S$ in cases (X1)--(X4) is 3-dimensional but not in (X6)--(X7).
Further, there exists an one-sided unit in~(X3) but not in~(X1),~(X2), and~(X4).
Since an automorphism of~$M_3(\mathbb{C})$ preserves a~rank of a matrix,
case~(X4) can not be isomorphic to neither~(X1) nor~(X2).
There exists an idempotent in~(X1) but not in~(X2)
such that its product with the second power of the radical of~$S$ is zero.
The radical of the subalgebra~$S_6$ from~(X6) is contained
in its own left annihilator from the whole~$S_6$;
on the other hand, the radical of the subalgebra~$S_7$ from~(X7) is not contained inside
its neither left nor right annihilator from~$S_7$.

\subsection{Case of unital 5-dimensional subalgebra. III}

Let $M = \Span\{e_{11},e_{12},e_{13},e_{22},e_{33}\}$.
An automorphism~$\varphi$~\eqref{AutPreserv} was actually constructed in Lemma~1
as the one preserving the radical $\Span\{e_{12},e_{13}\}$.
Thus, we may conclude that an automorphism $\varphi\in\Aut(M_3(\mathbb{C}))$
preserving~$M$ up to the action of $\Theta_{23}$ has the form~$\varphi$~\eqref{AutPreserv}
with $\lambda = \mu = 0$.

{\bf Theorem 7}.
Consider a direct decomposition of $M_3(\mathbb{C})$
with two subalgebras of the dimensions~4 and~5,
where a 5-dimensional subalgebra up to transpose and action of $\Aut(M_3(\mathbb{C}))$
is isomorphic to $M = \Span\{e_{11},e_{12},e_{13},e_{22},e_{33}\}$.
Then such decomposition up to transpose and up to action
of $\Aut(M_3(\mathbb{C}))$ is isomorphic to $S\oplus M$,
where $S$~is one of the following subalgebras:

(Y1) $S = \Span\{e_{21},e_{31},e_{32}+e_{33},e_{22}+e_{23}\}$;

(Y2) $S = \Span\{e_{21},e_{31},e_{22}+e_{32},e_{23}+e_{33}\}$;

(Y3) $S = \Span\{e_{21},e_{31},e_{11}+e_{22}+e_{23},e_{11}+e_{32}+e_{33}\}$;

(Y4) $S = \Span\{e_{21},e_{31},e_{11}+e_{22}+e_{32},e_{11}+e_{23}+e_{33}\}$;

(Y5) $S = \Span\{e_{21}+e_{22},e_{31}-e_{33},e_{32}+e_{33},e_{22}+e_{23}\}$;

(Y6) $S = \Span\{e_{21}-e_{31},e_{12}+e_{13}+e_{31},e_{11}+e_{32}+e_{33},e_{11}+e_{22}+e_{23}\}$;

(Y7) $S = \Span\{e_{21}-e_{23},e_{11}+e_{21}+e_{31},e_{12}+e_{22}+e_{32},e_{13}+e_{23}+e_{33}\}$;

(Y8) $S = \Span\{e_{21}+e_{22}+e_{33},e_{21}+e_{31},e_{11}+e_{12}-e_{13}+e_{22}+e_{32},e_{22}-e_{23}+e_{32}-e_{33}\}$;

(Y9) $S = \Span\{e_{21}+e_{22}-(x+1)e_{23},xe_{11}+e_{21}+e_{31},xe_{12}+e_{22}+e_{32},xe_{13}+e_{23}+e_{33}\}$;

(Y10) $S = \Span\{e_{21}+de_{22},e_{11}+e_{31}+d(e_{12}+e_{32}),e_{22}+e_{23},e_{12}+e_{13}+e_{32}+e_{33}\}$;

(Y11) $S = \Span\{e_{21}-e_{31}+e_{22}+e_{33},c(e_{12}+e_{13})+e_{21}+e_{22},e_{11}+e_{12}+e_{32}+e_{33},e_{11}-e_{13}+e_{22}+e_{23}\}$.

{\sc Proof}.
The subalgebra~$S$ is generated by matrices
$$
v_1 = \begin{pmatrix}
a & b & c \\
1 & d & 0 \\
0 & 0 & e \\
\end{pmatrix},\quad
v_2 = \begin{pmatrix}
g & h & i \\
0 & j & 0 \\
1 & 0 & k \\
\end{pmatrix},\quad
v_3 = \begin{pmatrix}
m & n & s \\
0 & p & 0 \\
0 & 1 & q \\
\end{pmatrix},\quad
v_4 = \begin{pmatrix}
u & v & x \\
0 & y & 1 \\
0 & 0 & z \\
\end{pmatrix}.
$$
Applying~$\varphi$~\eqref{AutPreserv} taken with
$\gamma = \lambda = \mu = 0$ and $\beta = \kappa a$, we may assume that $\underline{a = 0}$.

From the products
\begin{gather*}
v_1^2 = dv_1 + cv_4,\
v_1v_2 = gv_1 + ev_2 + iv_4, \
v_1v_3 = mv_1 + ev_3 + sv_4, \
v_1v_4 = uv_1 + (x+d)v_4, \\
v_2v_1 = jv_1 + bv_3, \
v_2^2 = (g+k)v_2 + hv_3, \
v_2v_3 = mv_2 + (k+n)v_3, \
v_2v_4 = uv_2 + vv_3 + jv_4, \\
v_3v_1 = pv_1 + v_2 + dv_3, \quad
v_3v_2 = qv_2 + jv_3, \quad
v_3^2 = (p+q)v_3, \quad
v_3v_4 = yv_3 + pv_4, \\
v_4v_1 = yv_1 + ev_4, \quad
v_4v_2 = v_1 + zv_2 + kv_4, \quad
v_4v_3 = zv_3 + qv_4, \quad
v_4^2 = (y+z)v_4,
\end{gather*}
we get the system of 69~equations,
\begin{gather*}
\{p,z\}\times\{q,y\} = 0,\quad qy+pz=1, \\
b = cu, \quad c(u-y) = 0, \quad e^2 - ed - cz = 0, \quad cv = 0, \quad c(e-d-x) = 0, \\
u(x+d) = 0, \quad xy + du - v = 0, \quad v(x+d) = 0, \quad e(z-u) - z(x+d) = 0,\\
cz - x(x+d) = 0, \quad
ey = 0,\quad v = eu, \quad v(d-e) = 0, \\
u(u-y-z) = 0, \quad v(u-z) = 0, \quad x(u-y) + v = 0, \\
em + sz = 0, \quad d(p-m) + n - ep - sy = 0, \\
b(p-m) + c - en -sv = 0, \quad c(q-m) - s(e+x) = 0, \\
j = -dp, \quad dm+g-n = 0, \quad e(q-p) - k - dq = 0, \quad c(m-p) + s(e-d) - i = 0, \\
m(m-p-q) = 0, \quad s(m-p) = 0, \quad n(m-q) + s = 0, \\
eg + iz = 0, \quad c = eg+iu, \quad h + d(j-g) - ej - iy = 0, \\
b(j-g) - eh - iv = 0, \quad c(k-g) - i(e+x) = 0, \\
bp = 0, \ e(k-j) + c - bq = 0, \ h = bm, \ b(g-j-n) + dh = 0, \ c(g-j) + ei - bs = 0, \\
hs = 0, \ h(j-k-n) = 0, \ i - gk - hq = 0, \ j(j-g-k) - hp = 0, \ h(m-q) = 0, \\
e = -kz, \quad c = i(u-z), \quad h(u-z) + v(j-k) - b = 0, \\
j(y-z) - yk - d = 0, \quad g(u-z) - uk + x = 0, \\
ju + pv = 0, \quad k(z-u) + x - vq - jz = 0, \\
v(g-j-n) + h(y-u) = 0, \quad x(g-j) + h - vs + i(z-u) = 0, \\
m(u-y) - pu = 0, \quad v(p-m) = 0, \quad x(m-p) + s(z-y) + n = 0, \\
nq + km - s = 0, \quad j(p-m) - p(k+n) = 0, \quad m(k+n) = 0, \\
n(g-k-n) + h(p-m) + i =0, \quad i(q-m) + s(g-k-n) = 0, \\
m(g-j) + s - qg = 0, \quad i(m-q) + s(k-j) = 0, \\
m(u-z) - uq = 0, \quad n(u-z) + v(p-q) + x = 0, \quad s(u-z) = 0.
\end{gather*}

{\sc Case I}: $p = z = 0$ and $y = 1/q$.
So, we have, $e = j = v = 0$, $m = uq^2$, $c = h$ and the remaining equations
(it can be verified by hands or with the help of Singular\cite{Singular})
\begin{gather*}
g = s/q, \quad x = -duq, \quad k = -dq, \quad u(u-1/q) = 0, \quad ug = 0, \\
c = iu, \quad b = cu, \quad n = dm+g, \quad i(m-q) + sk = 0, \quad dg - i(u-1/q) = 0.
\end{gather*}

{\sc Case IA}: $u,m \neq 0$.
Then $m = q$, $u = 1/q$, $g = s = 0$, $x = -d$, $n = dq = -k$, $i = cq$, and $b = c/q$.
It is the first solution required $q\neq0$
$$
v_1 = \begin{pmatrix}
0 & c/q & c \\
1 & d & 0 \\
0 & 0 & 0 \\
\end{pmatrix},\
v_2 = \begin{pmatrix}
0 & c & cq \\
0 & 0 & 0 \\
1 & 0 & -dq \\
\end{pmatrix},\
v_3 = \begin{pmatrix}
q & dq & 0 \\
0 & 0 & 0 \\
0 & 1 & q \\
\end{pmatrix},\
v_4 = \begin{pmatrix}
1/q & 0 & -d \\
0 & 1/q & 1 \\
0 & 0 & 0 \\
\end{pmatrix}.
$$
If $c = d = 0$, then we apply~$\varphi$~\eqref{AutPreserv} with
$\beta = \gamma = \lambda = \mu = 0$ and $\kappa = \nu q$ and get~(Y3).
If $d = 0$ and $c\neq0$, then we apply~$\varphi$ with
$\beta = \gamma = \lambda = \mu = 0$, $\kappa = \nu q$, and $\kappa\nu c = 1$ to get~(Y6).
If $d \neq 0$, then we apply~$\varphi$ with
$\beta = \gamma = \lambda = \mu = 0$, $\kappa = 1/d$, and $\nu = 1/(dq)$ to get~(Y11).

{\sc Case IB}: $u = m = 0$.
Then $b = c = x = 0$, $n = g = s/q$, $k = - dq$, and $i = -ds$.
It is the second solution required $q\neq0$
$$
v_1 = \begin{pmatrix}
0 & 0 & 0 \\
1 & d & 0 \\
0 & 0 & 0 \\
\end{pmatrix},\quad
v_2 = \begin{pmatrix}
s/q & 0 & -ds \\
0 & 0 & 0 \\
1 & 0 & -dq \\
\end{pmatrix},\quad
v_3 = \begin{pmatrix}
0 & s/q & s \\
0 & 0 & 0 \\
0 & 1 & q \\
\end{pmatrix},\quad
qv_4 = \begin{pmatrix}
0 & 0 & 0 \\
0 & 1 & q \\
0 & 0 & 0 \\
\end{pmatrix}.
$$
Let $s = 0$.
If $d = 0$, then we apply~$\varphi$ with
$\beta = \gamma = \lambda = \mu = 0$ and $\kappa = \nu q$ to get~(Y1).
If $d \neq 0$, then we apply~$\varphi$ with
$\beta = \gamma = \lambda = \mu = 0$, $\kappa = 1/d$, and $\nu = 1/(dq)$ to get~(Y5).

For $s\neq0$, we apply~$\varphi$ with
$\beta = \gamma = \lambda = \mu = 0$, $\nu = q/s$, and $\kappa = q^2/s$ to get~(Y10).

{\sc Case II}: $q = y = 0$ and $z = 1/p$.
Then $h = b = 0$, $m = up^2$,
\begin{gather*}
cu = 0, \quad v = du, \quad e - d - x + ku = 0, \quad
u(d-e) = 0, \quad u(u-1/p) = 0, \\
s = - emp, \quad g = j - k,\quad n = dm+g, \quad
i = - cp, \quad j = - dp, \quad k = - ep, \quad c = eg.
\end{gather*}

{\sc Case IIA}: $u,m\neq0$.
Then $m = p$, $u = 1/p$, $c = g = i = 0$, $e = d$,
$j = k = - dp = -n$, $v = d/p$, $x = -d$, and $s = - dp^2$.
It is the third solution with $p\neq0$
$$
v_1 = \begin{pmatrix}
0 & 0 & 0 \\
1 & d & 0 \\
0 & 0 & d \\
\end{pmatrix},\
v_2 = \begin{pmatrix}
0 & 0 & 0 \\
0 & -dp & 0 \\
1 & 0 & -dp \\
\end{pmatrix},\
v_3 = \begin{pmatrix}
p & dp & -dp^2 \\
0 & p & 0 \\
0 & 1 & 0 \\
\end{pmatrix},\
pv_4 = \begin{pmatrix}
1 & d & -dp \\
0 & 0 & p \\
0 & 0 & 1 \\
\end{pmatrix}.
$$
If $d = 0$, then we apply~$\varphi$ with
$\beta = \gamma = \lambda = \mu = 0$ and $p\nu = \kappa$ and get~(Y4).
If $d\neq0$, then we apply~$\varphi$ with
$\beta = \gamma = \lambda = \mu = 0$, $\kappa = 1/d$, and $p\nu = \kappa$ and get~(Y8).

{\sc Case IIB}: $u = m = 0$.
Then $s = v = 0$, $n = g$,
\begin{gather*}
x = e - d, \quad
i = - cp, \quad
j = - dp, \quad
k = - ep, \quad
g = p(e-d),\quad
c = eg.
\end{gather*}
It is the fourth solution
\begin{gather*}
v_1 = \begin{pmatrix}
0 & 0 & ep(e-d) \\
1 & d & 0 \\
0 & 0 & e \\
\end{pmatrix},\quad
v_2 = \begin{pmatrix}
p(e-d) & 0 & -ep^2(e-d) \\
0 & -dp & 0 \\
1 & 0 & -ep \\
\end{pmatrix},\\
v_3 = \begin{pmatrix}
0 & p(e-d) & 0 \\
0 & p & 0 \\
0 & 1 & 0 \\
\end{pmatrix},\quad
pv_4 = \begin{pmatrix}
0 & 0 & p(e-d) \\
0 & 0 & p \\
0 & 0 & 1 \\
\end{pmatrix}, \ p\neq0.
\end{gather*}
If $d = e = 0$, then we apply~$\varphi$ with
$\beta = \gamma = \lambda = \mu = 0$ and $\kappa = \nu p$ and get~(Y2).
If $d = 0$ and $e\neq0$, then we apply~$\varphi$ with
$\beta = \gamma = \lambda = \mu = 0$, $\kappa = 1/e$, $\nu = 1/(ep)$ to get~(Y7).
If $d \neq 0$, then we apply~$\varphi$ with
$\beta = \gamma = \lambda = \mu = 0$, $\kappa = 1/d$, and $\nu = 1/(dp)$ to get~(Y9).

Theorem is proved.\hfill $\square$

\subsection{Case of unital 5-dimensional subalgebra. IV}

Let $M = \Span\{e_{11}+e_{33},e_{12},e_{13},e_{22},e_{23}\}$.
As in the subsection~4.3, one can show that an automorphism
$\chi\in\Aut(M_3(\mathbb{C}))$ preserving~$M$
is defined by~\eqref{AutU}.

{\bf Theorem 8}.
Consider a direct decomposition of $M_3(\mathbb{C})$
with two subalgebras of the dimensions~4 and~5,
where a 5-dimensional subalgebra up to transpose and action of $\Aut(M_3(\mathbb{C}))$
is isomorphic to $M = \Span\{e_{11}+e_{33},e_{12},e_{13},e_{22},e_{23}\}$.
Then such decomposition up to transpose and up to action
of $\Aut(M_3(\mathbb{C}))$ is isomorphic to $S\oplus M$,
where $S$~is one of the following subalgebras:

(Z1) $S = \Span\{e_{11},e_{21},e_{31},e_{32}\}$;

(Z2) $S = \Span\{e_{11},e_{21},e_{31},e_{32}+e_{33}\}$;

(Z3) $S = \Span\{e_{11}+e_{22},e_{21},e_{31},e_{32}\}$;

(Z4) $S = \Span\{e_{11}+e_{22},e_{21}+e_{12},e_{31},e_{32}\}$.

{\sc Proof}. The subalgebra~$S$ is generated by matrices
$$
v_1 = \begin{pmatrix}
a & b & c \\
1 & d & e \\
0 & 0 & a \\
\end{pmatrix}, \quad
v_2 = \begin{pmatrix}
g & h & i \\
0 & j & k \\
1 & 0 & g \\
\end{pmatrix}, \quad
v_3 = \begin{pmatrix}
m & n & s \\
0 & p & q \\
0 & 1 & m \\
\end{pmatrix}, \quad
v_4 = \begin{pmatrix}
u & v & x \\
0 & y & z \\
0 & 0 & u+1 \\
\end{pmatrix}.
$$
Applying~$\psi$~\eqref{AutU} with $\beta = -\delta d$, $\varepsilon = \alpha p$, and
$\gamma = \varepsilon(a-d) -\alpha e$, we may assume that $\underline{d = e = p = 0}$.

From the products
\begin{gather*}
v_1^2 = av_1 - bv_4, \quad
v_1v_2 = gv_1 + av_2 - cv_4, \quad
v_1v_3 = mv_1 + av_3, \\
v_1v_4 = uv_1 + av_4, \
v_2v_1 = jv_1 + av_2 + bv_3 + (c-h)v_4, \
v_2^2 = kv_1 + 2gv_2 + hv_3, \\
v_2v_3 = mv_2 + (g+n)v_3 + sv_4, \quad
v_2v_4 = uv_2 + vv_3 + (g+x)v_4, \\
v_3v_1 = v_2 - nv_4, \quad
v_3v_2 = qv_1 + mv_2 + jv_3 + (k-s)v_4, \quad
v_3^2  = mv_3 + qv_4, \\
v_3v_4 = yv_3 + (z+m)v_4, \quad
v_4v_1 = yv_1 + (a-v)v_4, \quad
v_4v_2 = zv_1 + (u+1)v_2 + (g-x)v_4, \\
v_4v_3 = (u+1)v_3 + mv_4, \quad
v_4^2 = (2u+1)v_4,
\end{gather*}
we get the system of 80 equations, which radical computed by Singular~\cite{Singular}
gives $g = j = n = x = 0$ and the system
\begin{gather*}
\{b,m,u\}\times\{a,q,u+1\} = 0, \\
z(y+1) = 0, \quad
b(y+1) = 0, \quad
y(2u-y+1) = 0, \\
y(y^2-1) = 0, \quad
v(y-1) = 0, \quad
q(y-1) = 0, \\
z = my, \quad
v = ay, \quad
s = qv, \quad
k = aq, \quad
c = h = bm, \quad
i = ak + hm.
\end{gather*}
Because of the equation $u(u+1) = 0$, we may assume that $\underline{u = -1}$, otherwise
we apply $\Theta_{13}\circ T$.
Thus, $a = q = 0$, so, $k = s = v = 0$, and we have the system
\begin{gather*}
\{b,y\}\times \{y+1\} = 0,\quad
z = my, \quad
c = h = bm, \quad
i = bm^2.
\end{gather*}

{\sc Case I}: $y = -1$.
Then $z = -m$, and we have
$$
v_1 = \begin{pmatrix}
0 & b & bm \\
1 & 0 & 0 \\
0 & 0 & 0 \\
\end{pmatrix}, \
v_2 = \begin{pmatrix}
0 & bm & bm^2 \\
0 & 0 & 0 \\
1 & 0 & 0 \\
\end{pmatrix}, \
v_3 = \begin{pmatrix}
m & 0 & 0 \\
0 & 0 & 0 \\
0 & 1 & m \\
\end{pmatrix}, \
v_4 = -\begin{pmatrix}
1 & 0 & 0 \\
0 & 1 & m \\
0 & 0 & 0 \\
\end{pmatrix}.
$$
Applying $\psi$ with $\beta = \gamma = 0$ and $\varepsilon  = -\alpha m$,
we get the subalgebra
$S = \Span\{e_{21}+b'e_{12},e_{11}+e_{22},e_{31},e_{32}\}$, where $b' = \delta^2 b$.
If $b = 0$, then it is~(Z3).
Otherwise, we take~$\delta$ such that $\delta^2 b = 1$ to get~(Z4).

{\sc Case II}: $y = 0$.
Then $b = c = i = h = z = 0$, and
$S = \Span\{e_{11},e_{21},e_{31},e_{32}+me_{33}\}$.
If $m = 0$, then we get~(Z1).
If $m\neq0$, we apply~$\psi$ with $\beta = \gamma = \varepsilon = 0$
and $\delta = \alpha m$ to get~(Z2).

{\bf Remark 7}.
All cases (Z1)--(Z4) from Theorem~8  lie in different orbits under
action of automorphisms or antiautomorphisms of $M_3(\mathbb{C})$ preserving~$M$.
The radical in cases~(Z1) and~(Z3) is 3-dimensional but not in~(Z2) and~(Z4).
Further, there exists an one-sided unit in~(Z3) but not in~(Z1).
Finally, the radical of the subalgebra~$S_4$ from~(Z4) is contained
in its own left annihilator from the whole~$S_4$;
on the other hand, the radical of the subalgebra~$S_2$ from~(Z2) is not contained inside
its neither left nor right annihilator from~$S_2$.

\section*{Acknowledgments}

The work is supported by Mathematical Center in Akademgorodok.

\noindent Vsevolod Gubarev \\
Sobolev Institute of Mathematics \\
Acad. Koptyug ave. 4, 630090 Novosibirsk, Russia \\
Novosibirsk State University \\
Pirogova str. 2, 630090 Novosibirsk, Russia \\
e-mail: wsewolod89@gmail.com
\end{document}